\documentclass[11pt]{article}    % Specifies the document style.
\usepackage{amssymb,amsmath,amsthm,amscd}
\allowdisplaybreaks 

\begin{document}           % End of preamble and beginning of text.
\newtheorem{proposition}{Proposition}[section]
\newtheorem{theorem}[proposition]{Theorem}
\newtheorem{definition}{Definition}[section]
\newtheorem{corollary}[proposition]{Corollary}
\newtheorem{Lemma}[proposition]{Lemma}
\newtheorem{conjecture}{Conjecture}

\newtheorem{Remark}{Remark}[section]
\newtheorem{example}{Example}[section]

\newcommand{\ds}{\displaystyle}
\renewcommand{\ni}{\noindent}
\newcommand{\pa}{\partial}
\newcommand{\Om}{\Omega}
\newcommand{\om}{\omega}
\newcommand{\va}{\varepsilon}
\newcommand{\var}{\varphi_{y,\la}}

\newcommand{\la}{\lambda}
\newcommand{\sik}{\sum_{i=1}^k}
\newcommand{\vov}{\Vert\omega\Vert}
\newcommand{\Umy}{U_{\mu_i,y^i}}
\newcommand{\lamns}{\lambda_n^{^{\scriptstyle\sigma}}}
\newcommand{\chiomn}{\chi_{_{\Omega_n}}}
\newcommand{\ullim}{\underline{\lim}}
\newcommand{\bsy}{\boldsymbol}
\newcommand{\mvb}{\mathversion{bold}}
\newcommand{\R}{{\mathbb R}}
\newcommand{\bR}{{\mathbb R}}
\newcommand{\bC}{\mathbb{C}}
\newcommand{\bE}{\mathbb{E}}
\newcommand{\bH}{\mathbb{H}}
\newcommand{\bP}{\mathbb{P}}
\newcommand{\cF}{\mathcal{F}}

\newcommand{\beq}{\begin{eqnarray*}}
\newcommand{\eeq}{\end{eqnarray*}}

\newcommand{\ben}{\begin{enumerate}}
\newcommand{\een}{\end{enumerate}}

\newcommand{\beqs}{\begin{eqnarray*}&\displaystyle}
\newcommand{\eeqs}{&\end{eqnarray*}}

\renewcommand{\theequation}{\thesection.\arabic{equation}}
\catcode`@=11 \@addtoreset{equation}{section} \catcode`@=12
              % Produces the title.

%%%%%%%%%%%%%%%%%%%%%%%%%%%%%%%%%%%%%%%%%%%%%%%%%%%%%%%%%%%%%%%%%%%%

%%%%%%%%%%%%%%%%%%%%%%%%%%%%%%%%%%%%%%%%%%%%%%%%%%%%%%%%%%%%%%%%%%%%%

\title{SLE and $\alpha$-SLE driven by L\'evy processes\thanks{supported by EPSRC grant GR/T26368/01 in Britain and NSFC grant 10501048 in China.}}

\author{Qing-Yang Guan\thanks{Institute of Applied Mathematics, Academy of Mathematics and System Science, Chinese Academy of Sciences; email: guanqy@amt.ac.cn; presently: Department of Statistics, University of Oxford}\and Matthias Winkel\thanks{University of Oxford; email: winkel@stats.ox.ac.uk; supported by the Institute of Actuaries and Aon Limited}}

\maketitle

\begin{abstract} Stochastic Loewner Evolutions (SLE) with a multiple 
  $\sqrt{\kappa}B$ of Brownian motion $B$ as driving process are random planar 
  curves (if $\kappa\le 4$) or growing compact sets generated by a curve (if 
  $\kappa>4$). We consider here more general L\'evy processes as driving 
  processes and obtain evolutions expected to look like random trees or compact
  sets generated by trees, respectively. We show that when the driving force 
  is of the form $\sqrt{\kappa} B+\theta^{1/\alpha}S$ for a symmetric 
  $\alpha$-stable L\'evy process $S$, the cluster has zero or positive 
  Lebesgue measure according to whether $\kappa\le 4$ or $\kappa>4$. We also 
  give mathematical evidence that a further phase transition at $\alpha=1$ is 
  attributable to the recurrence/transience dychotomy of the driving L\'evy 
  process. We introduce a new class of evolutions that we call $\alpha$-SLE. 
  They have $\alpha$-self-similarity properties for $\alpha$-stable L\'evy 
  driving processes. We show the phase transition at a critical coefficient 
  $\theta=\theta_0(\alpha)$ analogous to the $\kappa=4$ phase transition.\\[12pt]
  \emph{AMS 2000 subject classifications: 60G51, 60G52, 60H10, 60J45.\newline
  }
\end{abstract}

\section{Introduction}

\em Loewner Evolutions \em are certain processes $(K_t)_{t\ge 0}$ taking values
in the space of closed bounded subsets of the complex upper half plane $\bH$ 
(or other simply connected domains), driven by a c\`adl\`ag function 
$U:[0,\infty)\rightarrow\bR$. They are best described via ordinary differential
equations
  \begin{equation} \partial_tg_t(z)=\frac{2}{g_t(z)-U(t)},\qquad g_0(z)=z,\qquad z\in\overline{\bH}=\{x+iy\in\bC:y\ge 0\},\label{leq}
  \end{equation}
as follows. $\partial_t$ is the right derivative as 
$U$ is right-continuous. For each $z\in\overline{\bH}$, the solution of 
(\ref{leq}) is well-defined on a time interval $[0,\zeta(z))$. Then the process
$K_t:=\{z\in\overline{\bH}:\zeta(z)\le t\}$, $t\ge 0$, is a strictly increasing
family of compact subsets of $\overline{\bH}$. We refer to $K_t$ as the \em 
cluster\em.

Loewner \cite{Loe-23} introduced these in the 1920s in a complex function 
theoretic framework of conformal mappings (the solutions 
$g_t:\bH\setminus K_t\rightarrow\bH$ of (\ref{leq}) are conformal mappings). In
the late 1990s, Schramm \cite{Sch-00} noticed that $U(t)=\sqrt{\kappa}B_t$ for 
a standard Brownian motion $B$ leads to an interesting class of \em Stochastic 
\em Loewner Evolutions $SLE_\kappa$, some of which he conjectured to be scaling
limits of important lattice models in statistical physics, subsequently proved
in collaboration with Lawler and Werner \cite{Law-01,Law-02} and by Smirnov
\cite{Smi-01}. Some introductory texts \cite{Law-04,Wer-02} are now available.
Cardy \cite{Car-03} gives a recent review of mathematical progress and further 
physical conjectures. 

Brownian motion is a suitable driving process since its independent identically
distributed (i.i.d.) increments translate into a composition of i.i.d.\ 
conformal mappings that describe, in a sense, independent growth increments.
Furthermore, Loewner evolutions transform well under Brownian scaling making
$SLE_\kappa$ conformally invariant, i.e.\ on the one hand, the
distribution of $(K_t)_{t\ge 0}$ is invariant under homotheties (the only 
conformal automorphisms of $\bH$ leaving start and end points $0$ and $\infty$ fixed), up to a 
linear time change; on the other hand, we can naturally consider $SLE_\kappa$
in other simply connected domains by application of a conformal mapping.

In this paper we discard the Brownian scaling property and
consider the larger class of processes with stationary independent increments
(L\'evy processes) as driving processes. Such processes are necessarily
discontinuous (except for Brownian motion, with drift). Whereas $SLE_\kappa$
is either a simple curve ($\kappa\le 4$) or generated by a curve
($\kappa>4$), \cite{Roh-01,Sch-00}, here, roughly, each discontinuity 
corresponds to a jump of the growth point on the boundary of the growing 
compact set. This leads to tree-like structures. Beliaev and Smirnov 
\cite{Bel-05} briefly mention such models in a complex analysis context as 
examples of fractal domains with high multifractal spectrum.

These models were recently introduced in the physics literature by Rushkin et 
al.\ \cite{RUS} who study driving processes of the form 
$U(t)=\sqrt{\kappa}B_t+\theta^{1/\alpha}S_t$ for a standard Brownian 
motion $B$ and an independent symmetric $\alpha$-stable L\'evy process $S$.
They observe two phase transitions.
\begin{enumerate}\item The Brownian phase transition of $SLE_\kappa$ at 
  $\kappa=4$ is not affected by the additional driving force 
  $\theta^{1/\alpha}S$. It can be expressed in terms of 
  $p(x)=\bP(\zeta(x)<\infty)$ as $p(x)=0$ for all $x\in\bR\setminus\{0\}$ for
  $\kappa\le 4$ versus $p(x)>0$ for all $x\in\bR\setminus\{0\}$ for $\kappa>4$.
  Due to the jumps, simulations look like trees and bushes respectively.
\item There is another phase transition at $\alpha=1$, which in the simulations
  yields ``isolated trees/bushes'' for $0<\alpha<1$ and ``forests of 
  trees/bushes'' for $1\le\alpha<2$.
\end{enumerate}
We strengthen their results from $x\in\bR$ to $z\in\overline{\bH}$ and 
rigorously establish the following theorem.

\begin{theorem}\label{thm1} Let $(K_t)_{t\ge 0}$ be an SLE driven by 
  $U_t=\sqrt{\kappa}B_t+\theta^{1/\alpha}S_t$ for a Brownian motion $B$ and 
  an independent symmetric $\alpha$-stable process $S$, with 
  $\zeta(z)=\inf\{t\ge 0:z\in K_t\}$. Then
  \begin{enumerate}\item[\rm (i)] if $0\le\kappa\le 4$ and $U\not\equiv 0$, 
    then for all 
    $z\in\overline{\bH}\setminus\{0\}$, we have $\bP(\zeta(z)=\infty)=1$;
  \item[\rm (ii)] if $\kappa>4$ and $1\le\alpha<2$, then for all 
    $z\in\overline{\bH}\setminus\{0\}$, we have $\bP(\zeta(z)<\infty)=1$;
  \item[\rm (iii)] if $\kappa>4$ and $0<\alpha<1$, then for all 
    $z\in\overline{\bH}\setminus\{0\}$, we have $0<\bP(\zeta(z)<\infty)<1$ 
    and $\lim_{z\rightarrow 0,z\in\overline{\bH}\setminus\{0\}}\bP(\zeta(z)<\infty)=1$.
  \end{enumerate}
\end{theorem}

Our methods combined with some probabilistic reasoning  
allow us to deduce the following corollary. Recall that L\'evy processes $C_t$ 
that are just the sums of finite numbers of jumps $\Delta C_s$ in any bounded 
interval $s\in[0,t]$ are called compound Poisson processes. A L\'evy process
$U$ is called recurrent if for all $a<0<b$ we have $\bE(\int_0^\infty 1_{\{a<U_t<b\}}dt)=\infty$, transient otherwise.

\begin{corollary}\label{cor1} Suppose that in the notation of the theorem, the driving 
  process is changed as follows, in terms of 
  $S_t^c=S_t-\sum_{s\le t}\Delta S_s1_{\{|\Delta S_s|>c\}}$, i.e. $S$ without 
  its big jumps, for some $c>0$, and independent compound Poisson processes 
  $R$ and $T$, recurrent and transient, respectively.
  \begin{enumerate}\item[\rm (i)] If $U_t=\sqrt{\kappa}B_t+\theta^{1/\alpha}S_t^c+R_t$ or $U_t=\sqrt{\kappa}B_t+\theta^{1/\alpha}S_t^c+T_t$, and $0\le\kappa\le 4$, but $\kappa>0$ or $\theta>0$ to avoid trivialities, then for all 
    $z\in\overline{\bH}\setminus\{0\}$, we have $\bP(\zeta(z)=\infty)=1$;
  \item[\rm (ii)] if $U_t=\sqrt{\kappa}B_t+\theta^{1/\alpha}S_t^c+R_t$ and $\kappa>4$ and $0<\alpha<2$, then for all 
    $z\in\overline{\bH}\setminus\{0\}$, we have $\bP(\zeta(z)<\infty)=1$;
  \item[\rm (iii)] if $U_t=\sqrt{\kappa}B_t+\theta^{1/\alpha}S_t^c+T_t$, and
    $\kappa>4$ and $0<\alpha<2$, then for all 
    $z\in\overline{\bH}\setminus\{0\}$, we have $0<\bP(\zeta(z)<\infty)<1$ 
    and $\lim_{z\rightarrow 0,z\in\overline{\bH}\setminus\{0\}}\bP(\zeta(z)<\infty)=1$.
  \end{enumerate}
\end{corollary}

This is strong evidence that the phase transition ``at $\alpha=1$'' is 
attributable to the recurrence/transience dychotomy of L\'evy 
processes. Under suitable regularity conditions on 
$\bP(|U_t|>x)\approx x^{-\alpha}$ as $x\rightarrow\infty$, such as regular 
variation, this is, of course, equivalent to $1\le\alpha\le\infty$ versus
$0<\alpha\le 1$, where a finer distinction is well-known at the critical 
value $\alpha=1$.  
%This Corollary strongly suggests the following conjecture.

%\begin{conjecture}\label{conj1} If $U_t$ is a L\'evy process with diffusion 
%  component $\sqrt{\kappa}B_t$ for some $\kappa\ge 0$, not piecewise 
%  deterministic, then
%  \begin{enumerate}\item[\rm (i)] if $0\le\kappa\le 4$, then for all 
%    $z\in\overline{\bH}\setminus\{0\}$, we have $\bP(\zeta(z)=\infty)=1$;
%  \item[\rm (ii)] if $\kappa>4$ and $\xi$ is recurrent, then for all 
%    $z\in\overline{\bH}\setminus\{0\}$, we have $\bP(\zeta(z)<\infty)=1$;
%  \item[\rm (iii)] if $\kappa>4$ and $\xi$ is transient, then for all 
%    $z\in\overline{\bH}\setminus\{0\}$, we have $0<\bP(\zeta(z)<\infty)<1$ 
%    and $\lim_{z\rightarrow 0,z\in\overline{\bH}\setminus\{0\}}\bP(\zeta(z)<\infty)=1$.
%  \end{enumerate}
%\end{conjecture}

Since recurrence and transience are governed only by rare big jumps, we \em 
expect \em that in the $\kappa\le 4$ case the phase transition is not reflected
in the local geometry of the cluster. Heuristically, in both cases pockets in 
the clusters will stabilise and remain unchanged after a while; in the 
transient case even the big trees themselves will remain unchanged eventually, 
whereas in the recurrent case bigger and bigger trees, possibly from the far 
left and the far right will almost meet above these unchanged pockets, and this
is reflected in the conformal mappings $g_t$ in that a whole pocket is mapped 
onto a very small portion of the upper half plane that ``disappears in the 
limit'' as $t\rightarrow\infty$; for $\kappa>4$ bigger bushes actually meet 
above pockets thereby incorporating the pockets in the cluster.

We leave the geometry of the cluster for further research, but
 establish the following result.

\begin{theorem}\label{thm1b} In the situation of Theorem \ref{thm1}, denote 
  Lebesgue measure on $\bH$ by $m$ and $B(0,r)=\{z\in\bH:|z|\le r\}$ for $r>0$.
  Then
  \begin{enumerate} \item[\rm (i)] if $0\le\kappa\le 4$, then 
    $m(\bigcup_{t\ge 0}K_t)=0$ a.s.;
    \item[\rm (ii)] if $\kappa>4$ and $1\le\alpha<2$, then 
      $m(\bH\setminus\bigcup_{t\ge 0}K_t)=0$ a.s.;
  \item[\rm (iii)] if $\kappa>4$ and $0<\alpha<1$, then 
    \beqs \lim_{r\downarrow 0}\frac{m\left(\bigcup_{t\ge 0}K_t\cap B(0,r)\right)}{m(B(0,r))}=1\qquad\mbox{and}\qquad \lim_{r\uparrow\infty}\frac{m\left(\bigcup_{t\ge 0}K_t\cap B(0,r)\right)}{m(B(0,r))}=0\qquad\mbox{a.s.}
    \eeqs
  \end{enumerate}
\end{theorem}
We actually believe that (ii) can be strengthened to $\bigcup_{t\ge 0}K_t=\overline{\bH}$
a.s. The other extreme is when the driving process is a compound Poisson 
process $U(t)=C_t$ with successive jump times $J_n$, $n\ge 1$, and jump heights
$X_n$, $n\ge 1$. $C$ is piecewise constant and hence the evolution can be 
decomposed and expressed as
\beqs g_{J_n+t}=\vartheta_{-X_1-\ldots-X_n}\circ g^0_t\circ\left(\vartheta_{X_n}\circ g^0_{J_n-J_{n-1}}\right)\circ\ldots\circ\left(\vartheta_{X_1}\circ g^0_{J_1}\right),\quad 0\le t<J_{n+1}-J_n,n\ge 0,
\eeqs
a composition of independent and identically distributed conformal mappings 
$\vartheta_{X_j}\circ g^0_{J_j-J_{j-1}}$, $j\ge 1$, where 
$g^0_t(z)=\sqrt{z^2+4t}$ is the conformal mapping from 
$\bH\setminus[0,2\sqrt{t}i]$ to $\bH$ that is associated with a driving 
function $U^0\equiv 0$ and $\vartheta_x(z)=z-x$ is a translation by 
$x\in\bR$. The flow $(\vartheta_{U_t}\circ g_t)_{t\ge 0}$ is similar to
flows of bridges (on $[0,1]$ instead of $\bH$) studied by Bertoin and Le Gall \cite{Ber-03}.

Clearly, $(K_t)_{t\ge 0}$ is here a forest of trees growing from $\bR$, with 
$g^0_{J_j-J_{j-1}}$ creating branches and $\vartheta_{X_j}$ moving the
growth point on the boundary.
Specifically, $K_t\cup\bR$ is path connected and, more precisely, has the \em
tree property \em that for all $y,z\in K_t\cup\bR$ there is a simple path 
$\rho:[0,1]\rightarrow\overline{\bH}$, unique up to time parameterisation, 
from $\rho(0)=y$ to $\rho(1)=z$ with $\rho(s)\in K_t\cup\bR$ for all 
$s\in[0,1]$. If $U$ is not a compound Poisson process, e.g. an $\alpha$-stable
L\'evy process, we have been unable to show that $K_t\cup\bR$ is path 
connected, but we believe, that the following holds.%\pagebreak

\begin{conjecture} If $U_t$ is a L\'evy process with diffusion component
  $\sqrt{\kappa}B_t$ for some $\kappa\ge 0$, then 
  \begin{enumerate}\item[\rm (i)] if $0\le\kappa\le 4$, then $K_t\cup\bR$ has
    the tree property for all $t\ge 0$. There is a simple left-continuous 
    function $\gamma:(0,\infty)\rightarrow\bH$ such that 
    $K_t\cap\bH=\{\gamma(s):0<s\le t\}$, for all $t\ge 0$.
  \item[\rm (ii)] if $\kappa>4$, then $K_t\cup\bR$ is generated by a 
    left-continuous function $\gamma:(0,\infty)\rightarrow\overline{\bH}$ in
    that $\bH\setminus K_t$ is the unbounded connected component of 
    $\bH\setminus\{\gamma(s):0<s\le t\}$, for all $t\ge 0$.
  \end{enumerate}
\end{conjecture}
This conjecture is a theorem for Brownian $SLE_\kappa$, see Rohde and
Schramm \cite{Roh-01} and Lawler et al. \cite{Law-01}, when $\gamma$ is indeed
continuous. In the setting of Theorem \ref{thm1}, the difficult part is to show
path connectedness of $\bR\cup K_t$, which is not obvious as the logarithmic 
spiral (see Marshal and Rohde \cite{Mar-05}) exemplifies. Heuristically, the
$\kappa=4$ phase transition is not affected by the small jumps since locally, 
the Brownian fluctuations dominate jump fluctuations as is expressed e.g. in 
$(U_{at}/\sqrt{a})_{t\ge 0}\rightarrow\sqrt{\kappa}B$ in distribution as $a\downarrow 0$, in the setting of the conjecture.

As a consequence of the scaling properties of (\ref{leq}) and Brownian motion
of the same index $2$, for $\theta=0$, any $\kappa\ge 0$ and $a>0$, the 
process $(\sqrt{a}K_t)_{t\ge 0}$, where $\sqrt{a}K_t=\{\sqrt{a}z:z\in K_t\}$, 
has the same distribution as $(K_{at})_{t\ge 0}$. The 
analogous statement for a pure $\alpha$-stable driving process, i.e.\ 
$\kappa=0$ and $\theta>0$ is not true: the distributions of 
$(a^{1/\alpha}K_t)_{t\ge 0}$ and $(K_{at})_{t\ge 0}$ are different. Scaling of 
index $2$ is intrinsic to equation (\ref{leq}).

However, we can construct clusters $(K_t)_{t\ge 0}$ such that 
$(a^{1/\alpha}K_t)_{t\ge 0}$ and $(K_{at})_{t\ge 0}$ have the same distribution
by modifying (\ref{leq}) to
\begin{equation} \partial_tg_t(z)=\frac{2|g_t(z)-U(t)|^{2-\alpha}}{g_t(z)-U(t)},\qquad g_0(z)=z,\qquad z\in\overline{\bH}=\{x+iy\in\bC:y\ge 0\},\label{leqalpha}
\end{equation}
for some $1<\alpha\le 2$. This equation still defines a process 
$(K_t)_{t\ge 0}$ of growing compact
subsets of $\overline{\bH}$, for a given c\`adl\`ag driving process $U$ and has
intrinsic scaling properties of index $\alpha$. We call this 
equation the $\alpha$-Loewner equation. The most interesting driving
processes are $\alpha$-stable processes, i.e.\ $\kappa=0$ in our setting. We
then derive the following phase transition.

\begin{theorem}\label{thm2} Let $1<\alpha<2$. If $(K_t)_{t\ge 0}$ is the $\alpha$-SLE driven by 
  $U_t=\theta^{1/\alpha}S_t$ for a symmetric $\alpha$-stable process $S$, then 
  there exists $\theta_0(\alpha)>0$ such that
  \begin{enumerate}\item[\rm (i)] if $0<\theta<\theta_0(\alpha)$, then for all 
    $z\in\overline{\bH}\setminus\{0\}$, we have $\bP(\zeta(z)=\infty)=1$;
  \item[\rm (ii)] if $\theta>\theta_0(\alpha)$, then for all 
    $z\in\overline{\bH}\setminus\{0\}$, we have $\bP(\zeta(z)<\infty)=1$.
  \end{enumerate}
\end{theorem}

Note that all driving processes are recurrent here, so the analogue to case
(iii) in the previous results does not arise. One could, however, e.g. add 
a transient compound Poisson process to the driving process and obtain the 
analogue to case (iii). We will also deduce the analogue of Theorem 
\ref{thm1b}.

\begin{corollary}\label{cor2} In the situation of Theorem \ref{thm2}, we have
  \begin{enumerate} \item[\rm (i)] if $0\le\theta<\theta_0(\alpha)$, then 
    $m(\bigcup_{t\ge 0}K_t)=0$ a.s.;
    \item[\rm (ii)] if $\theta>\theta_0(\alpha)$, then 
      $m(\bH\setminus\bigcup_{t\ge 0}K_t)=0$ a.s.
  \end{enumerate}
\end{corollary}

This class of growth processes $(K_t)_{t\ge 0}$ seems new and interesting. Theorem \ref{thm2}
and the discussion before describe some parallels to the class $SLE_\kappa$,
$\kappa\ge 0$. Our methods are strong enough to prove these analogous results,
even though the functions $g_t$ that solve (\ref{leqalpha}) are not conformal
mappings. The canonical driving processes are now jump processes, so we expect
the self-similar clusters to be trees or structures generated by trees. Again,
such structures are easily rigorously established for piecewise constant 
(e.g.\ compound Poisson) driving functions, but remain conjectural for stable
processes. It would be interesting to know if $\alpha$-SLE driven by 
$\alpha$-stable driving processes are scaling limits of natural lattice models.

The structure of this paper is as follows. In Section 2, we recall and extend
some preliminary results on fractional Laplacians, harmonic
functions and hitting time distributions; we also give an introduction to
Loewner evolutions and provide further and more detailed motivation for our 
class of driving functions. Sections 3 and 4 study the stochastic
differential equation of Bessel type that is associated with (\ref{leq}) for 
stochastic driving functions $U$ and deal with the proof of Theorem \ref{thm1} 
in the cases $z=x\in\bR$ and $z\in\bH$, respectively. In Section 5 we study the
increasing cluster $K_t$ and prove Theorem \ref{thm1b}. Section 6 is devoted to
properties of $\alpha$-SLE and the proof of Theorem \ref{thm2}.

\section{Preliminaries}

\subsection{Symmetric $\alpha$-stable processes and the fractional Laplacian}
\label{subsecsas}

Symmetric $\alpha$-stable L\'evy processes are Markov processes 
$(S_t)_{t\ge 0}$ starting from $S_0=0$, with \em stationary independent 
increments \em and c\`adl\`ag sample paths, whose distribution is given by
\beqs \bE(e^{i\lambda S_t})=e^{-t\psi(\lambda)},\qquad\psi(\lambda)=|\lambda|^\alpha=\int_{\bR\setminus\{0\}}(1-e^{i\lambda x}+i\lambda x1_{\{|x|\le 1\}})|x|^{-\alpha-1}dx
\eeqs
for some $0<\alpha<2$. We use Chapter VIII of Bertoin \cite{BER} as our
main reference. We can include $\alpha=2$, where $S_t=\sqrt{2}B_t$ is a 
Brownian motion $B_t$, and $S$ has as generator the Laplacian 
$\Delta_x=\partial_x^2$ on $\bR$. Brownian motion has the scaling property of 
index 2, called \em Brownian
scaling property \em that $(\sqrt{\kappa}B_t)_{t\ge 0}$ has the same 
distribution as $(B_{\kappa t})_{t\ge 0}$. For $0<\alpha<2$, the process $S$ 
has the \em scaling property of index $\alpha$ \em that 
$(\theta^{1/\alpha}S_t)_{t\ge 0}$ has the same distribution as 
$(S_{\theta t})_{t\ge 0}$. The infinitesimal generator of $S$ is the \em 
fractional Laplacian on $\R$\em, defined by the formula
\begin{align}\label{2.1}
\Delta^{\alpha/2}_x
w(x)=\lim_{\varepsilon\downarrow0}\mathcal{A}(1,-\alpha)\int_{\{x^\prime\in\bR:|x^\prime-x|>\varepsilon\}}\frac{w(x^\prime)-w(x)}{|x-x^\prime|^{1+\alpha}}\
dx^\prime,
\end{align}
where $w$ is a function on $\R$ such that the limit exists for all $x\in\R$, 
and $\mathcal{A}(1,-\alpha)$ is the constant 
$\alpha 2^{\alpha-1}\pi^{-1/2}\Gamma((1+\alpha)/2)/\Gamma(1-\alpha/2)$. We 
refer to Stein \cite{Ste-70} for an introduction and properties of the 
fractional Laplacian. We recall here that the domain of $\Delta^{\alpha/2}_x$ 
includes the Schwartz space of rapidly decreasing functions. It will be 
important in the sequel to apply (\ref{2.1}) as a \em formal generator \em
to functions where the limit does not exist for all $x\in\bR$, such as power
functions with a singularity at zero.

\begin{Lemma}\label{E:harmonic}
For $p\in\R$, define a function $w_p:\bR\rightarrow\bR$ by $w_p(0)=0$  and
 $$w_p(x)=|x|^{p-1},\  \  x\in\bR\setminus\{0\}, p\neq 1;\ \ \ \ w_1(x)=\ln|x|,\ \ x\in\bR\setminus\{0\}.$$
Then,
\begin{align}
\label{2.2}
&\Delta^{\alpha/2}_xw_p(x)=\mathcal{A}(1,-\alpha)\gamma(\alpha,p)
|x|^{p-\alpha-1},\ \ \ \mbox{for all $x\in\bR\setminus\{0\}$, and $p\in (0,\alpha+1)$,}
\end{align}
where $\gamma(\alpha,p)=\alpha^{-1}(p-1)\int_0^\infty
v^{p-2}(|v-1|^{\alpha-p}-(v+1)^{\alpha-p})\ dv$ for $p\neq1$ and
$\gamma(\alpha,1)=\alpha^{-1}\int_0^\infty
v^{-1}(|v-1|^{\alpha-1}-(v+1)^{\alpha-1})\ dv$.
\end{Lemma}\noindent{\bf Proof}$\ $
We assume without loss of generality that $x>0$. By definition
(\ref{2.1}) we have for $p\neq1$
\begin{align}\label{::1}
&\Delta^{\alpha/2}_xw_p(x)\nonumber\\
=&\lim_{\varepsilon\downarrow0}\mathcal{A}(1,-\alpha)\int_{\{x^\prime :|x^\prime -x|>\varepsilon\}}\frac{|x^\prime |^{p-1}-x^{p-1}}{|x-x^\prime |^{1+\alpha}}\
dx^\prime \nonumber\\
=&\lim_{\varepsilon\downarrow0}\mathcal{A}(1,-\alpha)
x^{p-\alpha-1}\int_{\{x^\prime :|x^\prime -1|>\varepsilon\}}\frac{|x^\prime |^{p-1}-1}{|x^\prime -1|^{1+\alpha}}\
dx^\prime \nonumber\\
=&\lim_{\varepsilon\downarrow0}\mathcal{A}(1,-\alpha)
x^{p-\alpha-1}\int_{\{x^\prime :|x^\prime |>\varepsilon\}}\frac{|x^\prime +1|^{p-1}-1}{|x^\prime |^{1+\alpha}}\
dx^\prime 
\nonumber\\
=&\lim_{\varepsilon\downarrow0}\mathcal{A}(1,-\alpha)
x^{p-\alpha-1}\int_\varepsilon^\infty\frac{|x^\prime +1|^{p-1}+|x^\prime -1|^{p-1}-2}{|x^\prime |^{1+\alpha}}\
dx^\prime 
\nonumber\\
=&\mathcal{A}(1,-\alpha)
\frac{(p-1)x^{p-\alpha-1}}{\alpha}\int_{\{x^\prime :x^\prime >0\}}\frac{(x^\prime +1)^{p-2}+(x^\prime -1)^{p-2}I_{\{x^\prime >1\}}-
(1-x^\prime )^{p-2}I_{\{0<x^\prime \leq1\}}}{|x^\prime |^{\alpha}}\
dx^\prime \nonumber\\
=&\mathcal{A}(1,-\alpha)
\frac{(p-1)x^{p-\alpha-1}}{\alpha}\int_0^\infty
v^{p-2}(|v-1|^{\alpha-p}-(v+1)^{\alpha-p})\ dv.
\end{align}
We use the transformation $(x^\prime +1)/x^\prime =v$ and $(x^\prime -1)/x^\prime =v$ in the last
step of (\ref{::1}). The case $p=1$ can be proved in the same way.\qed%\medskip
\pagebreak

\begin{Remark}\rm
\ By Lemma \ref{E:harmonic}, it is easy to check that $w_\alpha$ is
a harmonic function on $\R\setminus\{0\}$ for the symmetric
$\alpha$-stable process. When $\alpha>1$, $w_{\delta}$ is
subharmonic   and superharmonic  on $\R\setminus\{0\}$ when
$\delta\in(\alpha,\alpha+1)\cup(0,1)$ and $\delta\in[1,\alpha)$
respectively. When $0<\alpha<1$, $w_{\delta}$ is subharmonic
 and superharmonic on $\R\setminus\{0\}$ when
$\delta\in[1,\alpha+1)\cup(0,\alpha)$ and $\delta\in(\alpha,1)$
respectively. When $\alpha=1$, $w_{\delta}$ is a subharmonic
function on $\R\setminus\{0\}$ when $\delta\in(0,
1)\cup(1,\alpha+1)$.

By Lemma 4.2 in \cite{GQY1}, we can alternatively express the coefficients in 
Lemma 2.1 as
$\gamma(\alpha,p)=\int_0^1((u^{p-1}-1)(1-u^{\alpha-p})(1-u)^{-1-\alpha}
+(u^{p-1}-1)(1-u^{\alpha-p})(1+u)^{-1-\alpha})\ du$ for
$p\neq1$ and $\gamma(\alpha,1)=\int_0^1((1-u^{\alpha-1})\ln
(u)(1-u)^{-1-\alpha}+(1-u^{\alpha-1})\ln
(u)(1+u)^{-1-\alpha})\ du$. See also \cite[Lemma 5.1]{BBC}, 
\cite[Appendix]{IKE}, \cite[Appendix]{RUS} for other expressions of
these or closely related results.
% See also \cite{BBC} for a study of related (sub)harmonic functions on $(0,\infty)$.
  \end{Remark}

\subsection{Bessel-type processes and exit times}

Let $(B_t)_{t\geq0}$ and $(S_t)_{t\geq0}$ be standard Brownian
motion and an independent symmetric $\alpha$-stable process with generator
$\Delta^{\alpha/2}_x$, on a filtered probability space 
$(\Omega,\cF,(\cF_t)_{t\ge 0},\bP)$. Define
$U_t=\sqrt{\kappa}B_t+\theta^{1/\alpha}S_t$ and the conformal mappings
$(g_t)_{t\geq0}$ of SLE driven by $U_t$ via (\ref{leq}). 
%the following stochastic
%differential equation(see \cite{RUS}):
%$$\partial_tg_t(z)=\frac{2}{g_t(z)-U_t},\ \ \ g_0(z)=z,\ z\in
%\overline{\mathbb{H}}\setminus\{0\},$$where the derivative above is
%the right derivative as $U_t$ is right continuous. 
Let $h_t=g_t-U_t$, then we have the Bessel-type stochastic differential 
equation
\begin{align}\label{EQ:h_t} dh_t(z)=\frac{2dt}{h_t(z)} -dU_t,\ \ h_0(z)=z,\ z\in
\overline{\mathbb{H}}\setminus\{0\}.
\end{align}
$h_t(z)=h_{1,t}(z)+ih_{2,t}(z)$, $t\ge 0$, is an 
$\overline{\bH}$-valued Markov process, well-defined 
until hitting zero, for every $z\in\overline{\bH}\setminus\{0\}$ starting from
$z=z_1+iz_2$. The formal generator of the process $h$ is
\begin{align}\label{Plane:from cone::00}
Af(z)=\frac{-2z_2}{z_1^2+z_2^2}\partial_{z_2}f(z)+\frac{2z_1}{z_1^2+z_2^2}\partial_{z_1}f(z)+\frac{\kappa}{2}
\partial^2_{z_1}f(z)+ \theta\Delta^{\alpha/2}_{z_1}f(z).
\end{align}
It will be convenient to adopt a Markov process 
setup $(\Omega,\cF,(\cF_t)_{t\ge 0},(h_t)_{t\ge 0},%
(P_z)_{z\in\overline{\bH}\setminus\{0\}})$, slightly abusing notion, where 
$h_t$ under $P_z$ has the 
same distribution as $h_t(z)$ under $\bP$. In this vein, 
$\zeta=\inf \{t\geq 0: h_{t-}=0\mbox{ or }h_{t-}=U_t-U_{t-}\}$. We make a convention that $h_t=\Upsilon$,
a cemetery point $\Upsilon\not\in\overline{\bH}$, for $t\geq \zeta$ and $f(\Upsilon)=0$ for any function $f$.
 For a Borel set $D\subset\bH$,  denote $G_D(z,dz^\prime)=\int_0^\infty P^D_t(z,dz^\prime)dt$, where
$(P^D_t(z,dz^\prime))_{t\geq0}$ is the transition kernel for the  process
$(h_t)_{t\geq 0}$ killed when leaving $D$.

\begin{Lemma}\label{Hitting distribution}
Let $D$ be an open subset of $\mathbb{H}$ bounded away from $0$, i.e.\ such that $B(0,r)\subseteq
D^c$ for some $r>0$. Let $\tau=\inf\{t\geq 0: h_t
\notin D\}$ be the exit time from $D$, where $h_t$ is as in (\ref{EQ:h_t}). Then for every Borel set 
$B\subseteq  \overline{D}^c$ and every $z\in D$,
\begin{align}\label{Hitting distribution::0}
P_z\{h_{\tau}\in B\}=\ \int_DG_D(z,dz^\prime)\int_{\{z_1^{\prime\prime}\in\bR:z_1^{\prime\prime}+iz_2^\prime\in B\}}\frac{\theta\mathcal{A}(1,-\alpha)}{|z_1^{\prime\prime}-z_1^\prime|^{1+\alpha}}dz_1^{\prime\prime}.
\end{align}
where $z^\prime=z_1^\prime+iz_2^\prime$.
\end{Lemma}\noindent{\bf Proof}$\ $
We only need to prove that
\begin{align}\label{Hitting distribution::1}
E_zf(h_\tau)=\theta\mathcal{A}(1,-\alpha)\int_DG_D(z,dz^\prime)\int_{-\infty}^\infty
\frac{f(z_1^{\prime\prime}+iz_2^\prime)}{|z_1^{\prime\prime}-z_1^\prime|^{1+\alpha}}dz_1^{\prime\prime},
\end{align}
for each $C^2$ function $f$ on $\overline{\mathbb{H}}$ with compact
support satisfying
 $supp\ f\subseteq \overline{D}^c$. In fact by Dynkin's formula (see e.g.
It\^{o} \cite{Ito-05}),
 we have for all $z\in D$
\begin{align}\label{Hitting distribution::2}
E_zf(h_\tau)=&E_z\int_0^\tau Af(h_t)\ dt
=E_z\int_0^\tau\theta\Delta^{\alpha/2}_{z_1}f(h_t)\
dt\nonumber\\=&\int_0^\infty
\int_DP_t^D(z,dz^\prime)\theta\Delta^{\alpha/2}_{z_1^\prime}f(z^\prime)\
dt\nonumber\\=&\theta\mathcal{A}(1,-\alpha)\int_DG_D(z,dz^\prime)\int_{-\infty}^\infty
\frac{f(z_1^{\prime\prime}+iz_2^\prime)}{|z_1^{\prime\prime}-z_1^\prime|^{1+\alpha}}dz_1^{\prime\prime},\nonumber
\end{align}
which is (\ref{Hitting distribution::1}).\qed\medskip

Let $b>a>0$ and define ``inner'' and ``outer'' exit times of $h_{1,t}$ from 
$\{x\in\bR:a<|x|<b\}$ as
\begin{align}
\tau_{a,b}=\inf\{t\ge 0: |h_{1,t}|\leq a; |h_{1,s}|< b, \forall s\leq t \},\ \ \
\ \tau_{b,a}=\inf\{t\ge 0: |h_{1,t}|\geq b; |h_{1,s}|> a, \forall s\leq
t\},\end{align} where $\inf\varnothing=+\infty$. Let $\mu_{a,b}(z,dx^\prime)$
and $\mu_{b,a}(z,dx^\prime)$ be the conditional probability distributions 
under $P_z$ of $h_{1,\tau_{a,b}}$ and $h_{1,\tau_{b,a}}$ on events
$\{\tau_{a,b}<\infty\}$ and $\{\tau_{b,a}<\infty\}$ respectively.
Set $U_{a,b}=\{z\in\overline{\mathbb{H}}: a<\|z\|<b\}$,
where $\|z\|=\|z_1+iz_2\|=\max\{|z_1|,|z_2|\}$.  Denote similarly
\begin{align}
\overline{\tau}_{a,b}=\inf\{t\ge 0: \|h_t\|\leq a, \|h_s\|< b, \forall
s\leq t \},\ \ \ \ \overline{\tau}_{b,a}=\inf\{t\ge 0: \|h_t\|\geq b,
\|h_s\|> a, \forall s\leq t\},\end{align} 
and let
$\overline{\mu}_{a,b}(z,dx^\prime)$ and $\overline{\mu}_{b,a}(z,dx^\prime)$ be 
the conditional probability distributions of $h_{1,\overline{\tau}_{a,b}}$ and
$h_{1,\overline{\tau}_{b,a}}$ on events
$\{\overline{\tau}_{a,b}<\infty, h_{2,\overline{\tau}_{a,b}}\neq a\}$ and
$\{\overline{\tau}_{b,a}<\infty\}$ respectively.

\begin{Lemma}\label{ctu}
Let $b>a>0$, then the following assertions are true.
\begin{enumerate}\item[\rm(1)] Let $z\in
\mathbb{\overline{H}}$ such that $a<|z_1|<b$. Then $\mu_{a,b}(z,dx)$ may have
atoms at $x=a$ and $x=-a$ and is absolutely continuous on $\{x:|x|<a\}$ with 
density function $x\mapsto\varphi_{a,b}(z,x)$,
$\mu_{b,a}(z,dx)$ may have atoms at $x=b$ and $x=-b$, but is absolutely 
continuous on $\{x:|x|>b\}$ with density function $x\mapsto\varphi_{b,a}(z,x)$ 
such that
%\begin{align}
%I_{\{|x|<a/3\}}\mu_{a,b}(z,dx)=I_{\{|x|<a/3\}}\varphi_{a,b}(z,x)dx;\ \ \
%\ I_{\{|x|>2b\}}\mu_{b,a}(z,dx)=I_{\{|x|>2b\}}\varphi_{b,a}(z,x)dx,
%\end{align}
for all $|x|<a/3$, respectively $|x|>2b$,
\begin{align}\label{2.11}
\varphi_{a,b}(z,x)<\frac{3\cdot2^{3+4\alpha}}{a};\ \ \ \
\varphi_{b,a}(z,x)< 2^{3+4\alpha}\frac{(2b)^{\alpha}\alpha}{
|x|^{1+\alpha}}.
\end{align}

\item[\rm(2)] Let $z\in U_{a,b}\subset\overline{\bH}$.
Then $\overline{\mu}_{a,b}(z,dx)$ may have
atoms at $x=a$ and $x=-a$ and is absolutely continuous on $\{x:|x|<a\}$ with 
density function $x\mapsto\overline{\varphi}_{a,b}(z,x)$,
$\overline{\mu}_{b,a}(z,dx)$ may have atoms at $x=b$ and $x=-b$, but is 
absolutely continuous on $\{x:|x|>b\}$ with density function 
$x\mapsto\overline{\varphi}_{b,a}(z,x)$ 
such that
%Then there exist positive functions $\overline{\varphi}_{a,b}$,
%$\overline{\varphi}_{b,a}$ on $\R$ such that \begin{align}
%I_{\{|x|<a/3\}}\overline{\mu}_{a,b}(dx)=I_{\{|x|<a/3\}}\overline{\varphi}_{a,b}(x)dx;\
%\ \ \
%I_{\{|x|>2b\}}\overline{\mu}_{b,a}(dx)=I_{\{|x|>2b\}}\overline{\varphi}_{b,a}(x)dx,
%\end{align}
the same upper bounds as in (\ref{2.11}) hold.
\end{enumerate}

\end{Lemma}\noindent{\bf Proof}$\ $
 We only prove (2) as the proof of (1) is similar. Let
$|x|\geq|x'|\geq2b$. Then for any $|u|<b$, we have
\begin{align}\label{2.13}
2^{-2-2\alpha}\frac{|x'|^{1+\alpha}}{|x|^{1+\alpha}}\leq\frac{|x'-u|^{1+\alpha}}{|x-u|^{1+\alpha}}\leq
2^{2+2\alpha}\frac{|x'|^{1+\alpha}}{|x|^{1+\alpha}}.
\end{align}
Let $z\in \mathbb{\overline{H}}$ such that $z\in U_{a,b}$. For $|x|>b$, denote
\begin{align}\label{2.14}
f(x)=\frac{1}{P_z\{\overline{\tau}_{a,b}>\overline{\tau}_{b,a}\}} \int_{
U_{a,b}}\frac{\theta\mathcal{A}(1,-\alpha)}{|x-z_1^\prime|^{1+\alpha}}G_{U_{a,b}}(z,dz^\prime).
\end{align}
 By Lemma \ref{Hitting distribution}, we
know that $f$ is the density of $\overline{\mu}_{b,a}$ on
$\{x:|x|>b\}$. By (\ref{2.13}) and (\ref{2.14}), we see that for
$|x|>x'=2b$
\begin{align}\label{2.15}
2^{-2-2\alpha}\frac{(2b)^{1+\alpha}}{|x|^{1+\alpha}}f(2b)\leq
f(x)\leq 2^{2+2\alpha}\frac{(2b)^{1+\alpha}}{|x|^{1+\alpha}}f(2b).
\end{align}
Hence we have
$$2\int_{2b}^\infty 2^{-2-2\alpha}\frac{(2b)^{1+\alpha}}{|x|^{1+\alpha}}f(2b) dx\leq
\int_{-\infty}^{-2b}f(x) dx+ \int_{2b}^\infty f(x)dx\leq 1,$$ which
leads to $f(2b)\leq b^{-1}\alpha2^{2\alpha}$. Thus the
assertion concerning $\overline{\mu}_{b,a}$ follows from (\ref{2.15}).

Now let $|x|\leq|x'|\leq a/3$. Then for any $|u|>a$ we have
\begin{align}
2^{-2-2\alpha}\leq\frac{|u-x'|^{1+\alpha}}{|u-x|^{1+\alpha}}\leq
2^{2+2\alpha}.
\end{align}
Denote\begin{align}\label{ctu,,,}
f(x)=\frac{1}{P_z\{\overline{\tau}_{a,b}<\overline{\tau}_{b,a}, h_{2,\overline{\tau}_{a,b}}\neq a\}}
\int_{U_{a,b}}\frac{\theta\mathcal{A}(1,-\alpha)}{|z_1^\prime-x|^{1+\alpha}}G_{U_{a,b}}(z,dz^\prime),
\ \ \ |x|<a.
\end{align}
By definition of $\overline{\mu}_{a,b}$ and  Lemma \ref{Hitting
distribution}, we know that $f$ is the density of
$\overline{\mu}_{a,b}$ on $\{x:|x|<a\}$. By (2.16) and (2.17), we see
that for $|x|<x'=a/3$
\begin{align}\label{2.18}
2^{-2-2\alpha}f(a/3)\leq f(x)\leq 2^{2+2\alpha}f(a/3).
\end{align}
Hence we have
$$\int_{-a/3}^{a/3} 2^{-2-2\alpha}f(a/3)dx\leq
\int_{-a/3}^{a/3} f(x)dx\leq 1,$$ which leads to $f(a/3)\leq
3a^{-1}2^{1+2\alpha}$. Thus the assertion concerning $\overline{\mu}_{a,b}$ 
follows from (\ref{2.18}).
 \qed\medskip

\begin{Remark}\label{rfu}\rm
Let $g(x)=\ln|x|$ or $g(x)=|x|^{p-1}$ for $x\neq0$ and
$0<p<\alpha+1$. By Lemma \ref{ctu}, we see that $\int g \mu_{a,b}$,
$\int g\mu_{b,a}$, $\int g \overline{\mu}_{a,b}$ and $\int
g\overline{\mu}_{b,a}$ are all finite.

Whether conditional distributions such as $\mu_{a,b}$ have atoms at $a$ and
$-a$ depends on the so-called creeping properties of L\'evy processes (with 
drift), see Millar \cite{Mil-73} and Vigon \cite{Vig-04}. 
\end{Remark}

\subsection{Growing clusters, Loewner evolutions and independent increments}

The Riemann mapping theorem implies that for a 
compact set $K\subset\overline{\bH}$ such that $\bH\setminus K$ is simply 
connected, the family of conformal mappings $k:\bH\setminus K\rightarrow\bH$ 
is a set of three real dimensions. Since $\infty\not\in K$, it is natural to 
choose $k(\infty)=\infty$, the only point one can consistently fix for all 
compact sets $K$, with compositions of such conformal mappings in mind. The 
expansion at infinity then takes the form 
\beqs k(z)=a\left(z+b+\frac{{\rm hcap}(K)}{z}\right)+O\left(\frac{1}{z^2}\right),\qquad\mbox{for remaining parameters $a>0$ and $b\in\bR$,}
\eeqs
where ${\rm hcap}(K)$ is called the half-plane capacity (see Lawler 
\cite[Section 3.4]{Law-04}). It measures the size of $K$. Any increasing process $(K_t)_{t\ge 0}$ of compact sets with continuously 
increasing capacities can be (time-)parameterized such that 
${\rm hcap}(K_t)=2t$. Choosing $a=1$ is natural, $b=b_g:=0$ is one 
choice specifying a family of conformal mappings $(g_t)_{t\ge 0}$. Under the 
local growth condition 
\begin{equation} \bigcap_{\varepsilon>0}\overline{\left\{g_t(z):z\in K_{t+\varepsilon}\setminus K_t\right\}}=\left\{\mbox{single point}\right\}=:\left\{U(t)\right\}\qquad\mbox{for all $t\ge 0$},\label{locgrowth}
\end{equation}
where $\overline{C}$ denotes the closure of a Borel set 
$C\subset\overline{\bH}$, this growth point $b=b_h(t):=-U(t)$ is 
another choice for the parameter $b$ specifying another family of conformal
mappings $(h_t)_{t\ge 0}$. It can be checked that $(K_t)_{t\ge 0}$
is then the Loewner evolution driven by $(U(t))_{t\ge 0}$, the family 
$(g_t)_{t\ge 0}$ solves Loewner's differential equation (\ref{leq}), see 
Lawler \cite[Section 4.1]{Law-04}, and $h_t(z)=g_t(z)-U(t)$ solves the Bessel 
equation (\ref{EQ:h_t}) when integrating suitable test functions.
In general, $(U(t))_{t\ge 0}$ may be just measurable. However, we will assume 
in the sequel that $(U(t))_{t\ge 0}$ is c\`adl\`ag. The local growth condition,
even with a c\`adl\`ag function $(U(t))_{t\ge 0}$ is strictly weaker than the 
condition
\begin{equation} g_t^{-1}(\{U(t)\}):=\bigcap_{\varepsilon>0}\overline{g_t^{-1}(B(U(t),\varepsilon))}=\bigcap_{\varepsilon>0}\overline{K_{t+\varepsilon}\setminus K_t}=\{\mbox{single point}\}=:\{\gamma(t)\},\label{fngen}
\end{equation}
for a c\`adl\`ag function $\gamma:(0,\infty)\rightarrow\overline{\bH}$, where 
$B(x,\varepsilon)=\{z\in\bH:|z-x|\le\varepsilon\}$. In general, even 
under the local growth condition, equality may fail. If equality holds, one can
ask whether $(K_t)_{t\ge 0}$ is \em generated \em by a function $\gamma$ in a 
suitable class of functions, i.e. $\bH\setminus K_t$ is the unbounded connected
component of $\bH\setminus\overline{\{\gamma(s),0<s\le t\}}$, or even whether $\bH\cap K_t=\{\gamma(s-),0<s\le t\}$, i.e.
\begin{equation} \left\{z\in\bH\setminus K_{t-}:\lim_{\varepsilon\downarrow 0}g_{t-\varepsilon}(z)=U(t-)\right\}=\bH\cap K_t\setminus K_{t-}=\{\gamma(t-)\}.\label{swallow}
\end{equation}
In fact, $SLE_\kappa$ for 
$4<\kappa<8$ are examples where (\ref{fngen}) holds but (\ref{swallow}) 
fails -- further points in the left hand member of (\ref{swallow}) are 
called ``swallowed points''. The 
logarithmic spiral of Marshal and Rohde \cite{Mar-05} is an example where 
(\ref{fngen}) fails -- here the otherwise well-defined and continuous 
function $\gamma$ has neither left nor right limits at the time of the 
singularity, even though the driving function $(U(t))_{t\ge 0}$ is continuous. 
Werner \cite{Wer-02} remarks that one can build examples with a dense set of 
such singularities at different scales. In a rather more regular setting, it is
shown in \cite{Mar-05} that $1/2$-H\"older continuity of $(U(t))_{t\ge 0}$ with
small norm is sufficient for the existence and continuity of a simple curve 
$\gamma$. 
%We do not attempt to further investigate the existence and regularity of 
%$\gamma$ here.
% but do note the following result for the regularity of $U$.
%\begin{proposition} Under (\ref{locgrowth}), the function $t\mapsto U(t)$ is
%  always right-continuous with left limits.
%\end{proposition}
%\noindent{\bf Proof}$\ $ {\tt Is this really true? What if we start from 
%  $U(t)$ that is only measurable? If we can solve (\ref{leq}), why would 
%  $(K_t)_{t\ge 0}$ not satisfy the local growth condition???}
%\qed\medskip

Let us discuss further the geometric reasons for the choice of parameters, as
they provide further motivation for stochastic driving functions that are 
linear combinations of stable processes with stationary independent increments.
The first was $\infty\mapsto\infty$. Alternatively, one could fix $x\mapsto x$ 
for any specific $x\in\bR$, the boundary of $\overline{\bH}$, provided 
$x\not\in K$ but $K$ need not be compact. This is related to Loewner evolutions
``from $0$ to $x$'', rather than ``from $0$ to $\infty$''. 

Now let $(K_t)_{t\ge 0}$ be a Loewner evolution driven by any measurable 
function $(U(t))_{t\ge 0}$, growing ``from $0$ to $\infty$''; denote the 
associated solution to Loewner's equation by $(g_t)_{t\ge 0}$. The only 
conformal coordinate changes that leave zero and infinity fixed are homotheties
$z\mapsto cz$ inviting us to investigate 
$\widetilde{k}_t(z)=cg_t(z/c)$, $t\ge 0$. Clearly, these conformal mappings 
grow $(cK_t)_{t\ge 0}$, where ${\rm hcap}(cK_t)=c^2{\rm hcap}(K_t)$, so that we
reparameterise $k_t=\widetilde{k}_{c^{-2}t}$ and obtain
\begin{equation} \partial_tk_t(z)=\frac{2}{k_t(z)-cU_{c^{-2}t}},\qquad k_0(z)=z,\qquad z\in\overline{\bH},\label{leqscaling}
\end{equation}
so that $(cK_{c^{-2}t})_{t\ge 0}$ is a Loewner evolution driven by 
$(cU_{c^{-2}t})_{t\ge 0}$. This is the scaling property of index 2 that is
therefore intrinsic to Loewner's equation. 

\begin{proposition}[\cite{LSW-01,Roh-01} for $SLE_\kappa$]\label{propsii} 
  \begin{enumerate}\item[\rm(a)] An SLE $(K_t)_{t\ge 0}$ 
    is generated by a flow $h_t:\bH\setminus K_t\rightarrow\bH$ with stationary
    independent ``increments'' $h_{s,t}=h_t\circ h_s^{-1}$, $s\le t$, if and 
    only if the driving function $(U(t))_{t\ge 0}$ has the finite-dimensional 
    distributions of a L\'evy process.
  \item[\rm(b)] If $(U(t))_{t\ge 0}$ is a L\'evy process, then the distribution
    of $(\sqrt{a}K_{a^{-1}t})_{t\ge 0}$ is the same as that of $(K_t)_{t\ge 0}$
    if and only if $(U(t))_{t\ge 0}$ is a multiple of Brownian motion.
  \item[\rm(c)] If $U=\sqrt{\kappa}B+\theta^{1/\alpha}S$ for a Brownian motion $B$
    and an independent symmetric stable process of index $\alpha\in(0,2)$, then
    $(\sqrt{a}K_{a^{-1}t})_{t\ge 0}$ has the same distribution as a Loewner 
    evolution driven by 
    $\widetilde{U}=\sqrt{\kappa}B+\widetilde{\theta}^{1/\alpha}S$, where
    $\widetilde{\theta}=a^{\alpha/2-1}\theta$.
  \end{enumerate}
\end{proposition}
\noindent{\bf Proof}$\ $ For (a) just note that for fixed $s\ge 0$ and 
$h_t^{(s)}=h_{s+t}\circ h_s^{-1}$, we have by (\ref{EQ:h_t})
\begin{align}\nonumber dh_t^{(s)}(z)=dh_{s+t}(h^{-1}_s(z))=\frac{2dt}{h_{s+t}(h^{-1}_s(z))} -dU_{s+t}=\frac{2dt}{h_t^{(s)}(z)}-dU_t^{(s)},\ \ h^{(s)}_0(z)=z,\ z\in
\overline{\mathbb{H}}\setminus\{0\},
\end{align}
where $U_t^{(s)}=U_{s+t}-U_s$, and this easily yields the result. (b) 
and (c) are simple consequences of the scaling properties of Loewner's 
equation, (\ref{leqscaling}), and of $B$ and $S$ (see Subsection \ref{subsecsas}).
\qed\medskip

The property in (b) is called conformal invariance. For any simply connected 
domain $D\subset\bC$, $D\neq\bC$, one can now uniquely define $SLE_\kappa$ 
from one boundary point $\alpha$ to another boundary point $\beta$ by conformal mappings $f:\bH\rightarrow D$ with 
$f(0)=\alpha$ and $f(\infty)=\beta$, up to a 
linear time change. For any other L\'evy process, the 
definition is not unique. However, note that for the driving processes in (c), 
the properties of SLE studied in this paper do not depend on $\theta$.
%\pagebreak

\section{$\bR$-valued Bessel-type processes driven by $U=\sqrt{\kappa}B+\theta^{1/\alpha}S$}

By (\ref{EQ:h_t}), it is easy to see that $(h_t(x))_{0\le t<\zeta(x)}$ is
$\R$-valued for all $x\in\R\setminus\{0\}$. In this case their formal generator
$A$ reduces to
$$Af(x)=\frac{2}{x}\partial_xf(x)+\frac{\kappa}{2}\partial_x^2f(x)+\theta\Delta^{\alpha/2}_xf(x),\ \ \ \mbox{for all  $x\in\bR\setminus\{0\}$.}$$

\begin{proposition}\label{hitting, lemma }
\label{Line:<4}When $0\leq\kappa\leq 4$ and $0<\alpha<2$, we have
$\zeta(x)=\infty$ a.s.\ for all $x\in \R\setminus \{0\}$.
\end{proposition}\noindent{\bf Proof}$\ $
We will use the same notation as in Lemma \ref{E:harmonic} and always
assume that $\kappa>0$. The case $\kappa=0$ can be proved similarly. 
%Without loss of generality we assume $x>0$.

\noindent Case 1. $0<\alpha\leq1$. By Lemma \ref{E:harmonic}, we
have for $y\in\bR\setminus\{0\}$
$$Aw_1(y)=\frac{2}{y}\partial_yw_1(y)+\frac{\kappa}{2}\partial^2_yw_1(y)
+\theta\Delta^{\alpha/2}_yw_1(y)\geq\theta\Delta^{\alpha/2}_y
w_1(y)=\theta\mathcal{A}(1,-\alpha)\gamma(\alpha,1)
|y|^{-\alpha}\geq 0.$$ For $0<a<b$, let $ \tau_{a,b}$ and
$\tau_{b,a}$ be the inner and outer exit times defined in (2.8). Let 
$\mu_{a,b}$ and
$\mu_{b,a}$ be the corresponding conditional probability
distribution.
 By Dynkin's formula we have
\begin{align}
&\ln |x|\leq P_x\{\tau_{a,b}<\tau_{b,a}\}\int_{\{|y|\leq a\} }\ln |y|\
{\mu_{ a,b}(x,dy)}+P_x\{\tau_{a,b}>\tau_{b,a}\}\int_{\{|y|\geq b \}}\ln
|y|\ {\mu_{ b,a}(x,dy)}.\nonumber
\end{align}
Therefore
\begin{align}\label{111}
 P_x\{\tau_{a,b}<\tau_{b,a}\}\leq
\frac{\ln|x|-\int_{\{|y|\geq b\} }\ln |y|\ {\mu_{
b,a}(x,dy)}}{\int_{\{|y|\leq a\}  }\ln |y| {\mu_{a,b}(x,dy)}-\int_{\{|y|\geq b\}
}\ln |y|\ {\mu_{ b,a}(x,dy)}}.
\end{align}
By Lemma \ref{ctu} we know that $\int_{\{|y|\geq b\}}\ln |y|\ {\mu_{
b,a}(x,dy)}$ is bounded for fixed $b$ uniformly in $a<b$. Letting
$a\downarrow0$ in (\ref{111}) we get $\zeta=\infty$, $P_x$-a.s.\medskip

\noindent Case 2. $0<\kappa<4, 1<\alpha<2 $. Let
$f_1=w_{3/2-2/\kappa}$. First we prove the case
$\kappa\geq 2$. By Lemma \ref{E:harmonic} we have for $y\neq 0$
\begin{align}
Af_1(y)&=\left(\frac{2}{y}\partial_y
+\frac{\kappa}{2}\partial^2_y\right)w_{3/2-2/\kappa}(y)+\theta\Delta^{\alpha/2}_yw_{3/2-2/\kappa}
(y)
\nonumber\\&=\left(\frac{1}{2}-\frac{2}{\kappa}\right)\left(1-\frac{\kappa}{4}\right)|y|^{-3/2-2/\kappa}+
\theta\mathcal{A}(1,-\alpha)\gamma\left(\alpha,\frac{3}{2}-\frac{2}{\kappa}\right)
|y|^{1/2-2/\kappa-\alpha}.
\end{align}
Noticing that $(\frac{1}{2}-\frac{2}{\kappa})(1-\frac{\kappa}{4})<0$
we can find a constant $c$ such that $Af_1(y)-cf_1(y)<0$ for all $y\neq
0$. Again by Dynkin's formula we obtain
\begin{align}\label{Line:<4:1}
f_1(x)&\geq E_x\left[e^{-c\tau_{a,b}}f_1(h_{\tau_{a,b}})\right]+
E_x\left[e^{-c\tau_{b,a}}f_1(h_{\tau_{b,a}})\right].
\end{align}
If $P_x\{\zeta<\infty\}>0$, we can choose   $b,T\in\R$ big enough
such that $P_x\{\lim_{a\downarrow0}\tau_{a,b}<T\}>0$. Hence by
(\ref{Line:<4:1}), we get $f_1(x)\geq
e^{-cT}P_x\{\lim_{a\downarrow0}\tau_{a,b}<T\}a^{1/2-2/\kappa}+
E_x[e^{-c\tau_{b,a}}f_1(h_{\tau_{b,a}})]$, which is impossible when
taking $a\downarrow0$. When $0<\kappa<2$, we can take
$f_1=w_{\frac{1}{2}}$ and use the same method.\medskip

\noindent Case 3. $\kappa=4, 1<\alpha<2 $. By Lemma \ref{E:harmonic}
we have $(\frac{2}{y}\partial_y +2\partial^2_y)w_{1}(y)=0$. Therefore for $y\neq 0$
and $c>0$ we have
 \begin{align}\label{Line:<4:2}
A(w_1+cw_{3-\alpha})(y)%\nonumber\\
=&c\left(\frac{2}{y}\partial_y
+2\partial^2_y\right)w_{3-\alpha}(y)+\theta\Delta^{\alpha/2}_yw_1(y)+c\theta\Delta^{\alpha/2}_yw_{3-\alpha}(y)
\nonumber\\
=&2c(2-\alpha)^2|y|^{-\alpha}+
\theta\mathcal{A}(1,-\alpha)\gamma(\alpha,1) |y|^{-\alpha}+
c\theta\mathcal{A}(1,-\alpha)\gamma(\alpha,3-\alpha)
|y|^{2-2\alpha}.
\end{align}
By (\ref{Line:<4:2})  and noticing that $-\alpha<2-2\alpha$, we can
find $c$ large enough  and  $r>0$ small enough such that $Af_2(y)>0$
for $|y|<r$, $y\neq 0$. Then  following the same method as in case 1 we
can prove  $P_x\{\tau_{0,r}<\tau_{r,0}\}=0$, which leads to the
conclusion.\qed

\begin{proposition}\label{Line:>4}
When $4<\kappa$ and $1\leq\alpha<2$, we have $\zeta(x)<\infty$ a.s.
for all $x\in \R\setminus \{0\}$.
\end{proposition}\noindent{\bf Proof}$\ $
We will use the same notation in Lemmas \ref{E:harmonic} and
2.3. Without loss of generality we assume $x>0$.

\noindent Case 1. $2-4/\kappa\leq\alpha<2$. In this case
$\gamma(\alpha, 2-4/\kappa)\le 0$. We get by Lemma
\ref{E:harmonic} that $Aw_{2-4/\kappa}\leq 0$. By Dynkin's
formula we have
\begin{align}
P_x\{\tau_{a,b}<\tau_{b,a}\}\geq \frac{\int_{\{|y|\geq b\}
}|y|^{1-4/\kappa}\ {\mu_{
b,a}(x,dy)}-|x|^{1-4/\kappa}}{\int_{\{|y|\geq b\} }
|y|^{1-4/\kappa}\ {\mu_{ b,a}(x,dy)}-\int_{\{|y|\leq a\}  }
|y|^{1-4/\kappa} {\mu_{ a,b}(x,dy)}}.
\end{align}
By Lemma \ref{ctu}, letting $a\downarrow 0$ and then $b\uparrow \infty$ we
get the conclusion.\medskip

\noindent Case 2. $1<\alpha<2-4/\kappa$. By Lemma
\ref{E:harmonic}, we can check $Aw_{\alpha}<0$. Hence we  can get
the same conclusion by the method above.\medskip

\noindent Case 3. $\alpha=1$. By Lemma \ref{E:harmonic}, we can
check that there exists a number $c>0$ satisfying
$Aw_{3/2-2/\kappa}(y)<0$ for $0<|y|<c$. Hence we
obtain $\lim_{y\downarrow 0}P_y\{\tau_{0,c}<\tau_{c,0}\}=1$ by Dynkin's
formula.
 Now, by the Markov property, we only need to prove that $P_x\{\tau_{a,\infty}<\infty\}=1$ for all $a>0$ and  $x\neq 0$.
Here $\tau_{a,\infty}=\inf_{b>a}\tau_{a,b}$.

By Lemma \ref{E:harmonic}, we have $Aw_1(y)<0$ for $y\neq 0$. Hence
we have by Dynkin's formula
\begin{align}
 P_x\{\tau_{a,b}<\tau_{b,a}\}\geq
\frac{\ln|x|-\int_{\{|y|\geq b \}}\ln |y|\ {\mu_{
b,a}(x,dy)}}{\int_{\{|y|\leq a \} }\ln |y| {\mu_{a,b}(x,dy)}-\int_{\{|y|\geq b\}
}\ln |y|\ {\mu_{ b,a}(x,dy)}}.\nonumber
\end{align}
By Lemma 2.3, letting $b\uparrow \infty$ we have
$P_x\{\tau_{a,\infty}<\infty\}=1$. \qed

\begin{Lemma}\label{ES:Line:<4:<1}
Let $4<\kappa$ and $0<\alpha<1$. There exist constants $k_1,k_2>0$
depending on $\kappa, \alpha, \theta$ such that
\begin{align}P_x\{\zeta=\infty\}>k_2,\qquad\mbox{for all $x\geq k_1$.}\end{align}
\end{Lemma}\noindent{\bf Proof}$\ $
By Lemma \ref{E:harmonic}, we can choose $c$ large enough such that
$Aw_{\alpha/2+1/2}(y)<0$ for $|y|>c/2$. Hence we have
\begin{align}
P_x\{\tau_{c/2,b}>\tau_{b,c/2}\}\geq
\frac{\int_{\{|y|\leq c/2\}}|y|^{\alpha-1}\ {\mu_{
b,a}(x,dy)}-c^{\alpha-1}}{\int_{\{|y|\leq c/2\}} |y|^{\alpha-1}
{\mu_{ a,b}(x,dy)}-\int_{\{|y|\geq b\} } |y|^{\alpha-1}\ {\mu_{
b,a}(x,dy)}},\ \ \ c<x<b.\nonumber
\end{align}
By Lemma 2.3, letting $b\uparrow \infty$, we get the conclusion.
\qed

\begin{proposition}\label{Line:<4:<1}
Let $4<\kappa$ and $0<\alpha<1$. There exists constant $c>0$ such
that
\begin{align}
\noindent &(\mathrm{a})\ \
\frac{1}{c}|x|^{1-4/\kappa}<P_x\{\zeta=\infty\}<c|x|^{1-4/\kappa},\
\ 0<|x|\leq 1;\nonumber\\ &(\mathrm{b})\ \
\frac{1}{c}|x|^{\alpha-1}<P_x\{\zeta<\infty\}<c|x|^{\alpha-1},\ \
|x|>1.\nonumber\end{align}
\end{proposition}
\noindent{\bf Proof}$\ $ First we prove the upper bound in (a). Define 
functions $u_1(y)=|y|^{1-{2}/{\kappa}}\wedge 2$ and
$u_2(y)=|y|^{1-{4}/{\kappa}}\wedge 2$. Now we suppose
$1-2/\kappa<\alpha$. By Lemma \ref{E:harmonic} and direct
calculation we have
\begin{align}\label{Line:<4:<1::1}
\frac{1}{c_1} <\lim_{|y|\downarrow
0}\Delta^{\alpha/2}_yu_1(y)/|y|^{1-{2}/{\kappa}-\alpha}<c_1;\ \ \ \
\ \frac{1}{c_2} <\lim_{|y|\downarrow
0}\Delta^{\alpha/2}_yu_2(y)/|y|^{1-{4}/{\kappa}-\alpha}<c_2.
\end{align}
for some positive constants $c_1$ and $c_2$. Choose a small positive real
number $c_3$ such that $u_2(y)-c_3u_1(y)>0$ for $y\neq0$. By Lemma
\ref{E:harmonic} we have
\begin{align}\label{Line:<4:<1::2}
 A(u_2- c_3u_1)(y)&=-c_3(1-{2}/{\kappa})|y|^{-1-{2}/{\kappa}}
+\theta\Delta^{\alpha/2}_y(u_2-c_3u_1)(y).
\end{align}
Let $f_1=w_{2-{4}/{\kappa}}- c_3u$. By (\ref{Line:<4:<1::1}) and
(\ref{Line:<4:<1::2}), we can find a positive real number $c_4$ such
that $Af_1(y)<0$ for $y\neq0$ and $|y|<c_4$. Applying the same
notation as in Proposition \ref{Line:<4}, we have for $0<a<c_4$
\begin{align} P_x\{\tau_{a,c_4}>\tau_{c_4,a}\}\leq
\frac{f_1(x)}{\int_{|y|\geq c_4 }f_1(y)\ {\mu_{
c_4,a}(x,dy)}-\int_{|y|\leq a  }f_1(y) {\mu_{a,c_4}(x,dy)}}.\nonumber
\end{align}
By Lemma 2.3, letting $a\downarrow 0$ in the equality above, we have
\begin{align}
  P_x\{\zeta=\infty\}\leq P_x\{\tau_{0,c_4}>\tau_{c_4,0}\}\leq
\frac{x^{1-{4}/{\kappa}}}{\lim_{a\downarrow 0}  \int_{|y|\geq c_4
}f_1(y)\ {\mu_{ c_4,a}(x,dy)}},\nonumber
\end{align}which gives the second inequality in (a). When $1-2/\kappa\geq
\alpha$, we can prove the upper bound in the same way as above by
noticing that
\begin{align}
\frac{1}{c} <\lim_{|y|\downarrow
0}\Delta^{\alpha/2}_yu(y)/\ln|y|<c,\qquad\mbox{when $\beta= \alpha$;}\nonumber\\
|\Delta^{\alpha/2}_yu(y)|<c,\quad y\in (-1,1),\qquad\mbox{when $\beta>\alpha$,}
\end{align}
for some constant $c$ depending on $\beta$ and $\alpha$, where
$u(y)=|y|^\beta\wedge2$. This can be checked directly, see also
Proposition 2.3 in \cite{GQM} and Proposition 2.5 in \cite{GQY1}.

Next we prove the lower bound in (a). We use the notation $k_1$ and $k_2$
as in Lemma \ref{ES:Line:<4:<1}. Let
$u_3(y)=|y|^{1-4/\kappa}\wedge M$ for some $M>0$. Choose $M$ big
enough such that $Au_3(y)>0$ for $0<|y|<k_1$.  By this fact and
applying the same method as above, we can prove that for some
constant $c_5$
\begin{align}
P_x\{\tau_{k_1,0}<\tau_{0,k_1}\}\geq c_5|x|^{1-{4}/{\kappa}},\qquad
0<x<k_1.\nonumber
\end{align}  
Hence by the Markov property and Lemma \ref{ES:Line:<4:<1} we get
$P_x\{\zeta=\infty\}\geq k_2c_5|x|^{1-{4}/{\kappa}}$ and complete
the proof of (a). We omit the proof of (b) as it can be proved  by
similar discussions.\qed\medskip

\section{$\overline{\bH}$-valued Bessel-type processes driven by $U=\sqrt{\kappa}B+\theta^{1/\alpha}S$}

In this section we consider the problem whether the Bessel-type process on the 
complex upper half plane, given in (\ref{EQ:h_t}), can hit $0$. Denote this
process by $h_t(z)=h_{1,t}(z)+ih_{2,t}(z)$ and $z=z_1+iz_2$. For
$z\in \mathbb{\overline{H}}$, we have that
\begin{align}\label{EQ:h_1:h_2}
\left\{ \displaystyle \begin{array}{ll}
   \displaystyle dh_{1,t}(z)=&\displaystyle\frac{2h_{1,t}(z)dt}{h_{1,t}^2(z)+h_{2,t}^2(z)} -dU_t,\ \
h_{1,0}(z)=z_1,
   \\
 \ds   dh_{2,t}(z)=&\displaystyle\frac{-2h_{2,t}(z)dt}{h_{1,t}^2(z)+h_{2,t}^2(z)},\ \ \  \ \ \ \ \ \ \
h_{2,0}(z)=z_2.\ \ \ \ \ \
\end{array}\right. \end{align}

\subsection{The subcritical phase $0<\kappa<4$}

We have to prepare some results to deal with the hitting problem. For
$\delta>0$, denote by $V_\delta=\{z=z_1+iz_2:0<z_2\leq \delta|z_1|)\}$ the
double wedge of slope $\delta$,
and $\tau_\delta=\inf\{t\geq 0: h_t\in V_\delta\}$ the first entrance time.

\begin{Lemma}\label{Plane:hitting cone}
If $\kappa>0$, then for each $\delta>0$ and $z\in{\mathbb{H}}$,
\begin{align}\label{Plane:hitting cone::1}
P_z\{\tau_\delta<\infty\}=1.
\end{align}\end{Lemma}
\noindent{\bf Proof}$\ $ The proof is in five parts.\medskip

\noindent \em 1. We reduce the proof to small $z$. \em We only need to prove 
(\ref{Plane:hitting cone::1}) when $z\notin V_\delta$. Without loss of 
generality we
assume that $\delta<1$. Let $s>0$ and  denote
$${d}_{\delta,s}=\inf\{t\geq0: h_t \in V_\delta\mbox{ or }h_{2,t}\leq s \}.$$
We claim that $d_{\delta,s}<\infty$ a.s. This will follow if we show
$P_z(E)=0$ for
$$\ds E=\left\{\omega\in\Omega: |h_{1,t}(\omega)|<z_2/\delta, \ h_{2,t}(\omega)>s\mbox{ for all $t>0$}\right\}. $$
 In fact, we have for a.e. $\omega\in E$
  $$\ds\lim_{t\rightarrow\infty}h_{2,t}(\omega)=z_2+\lim_{t\rightarrow\infty}\int_0^t
  \frac{-2h_{2,t}(\omega)dt}{h_{1,t}^2(\omega)+h_{2,t}^2(\omega)}\leq z_2
  -\ds\lim_{t\rightarrow\infty}\int_0^t\frac{2sdt}{z_1^2/\delta^2+z_2^2}=-\infty,$$
  which is absurd for a process in $\overline{\bH}$. Next, by the Markov 
  property,
\begin{align}\label{Plane:hitting cone::2}
P_z\{\tau_\delta<\infty\}=&P_z\{h_{d_{\delta,s}}\in
V_\delta\}+P_z\{h_{d_{\delta,s}}\notin V_\delta,
\tau_\delta<\infty\}\nonumber\\
=&P_z\{h_{d_{\delta,s}}\in V_\delta\}+E_z\left[I_{\{h_{d_{\delta,s}}\notin
V_\delta\}} P_{h_{d_{\delta,s}}}\{\tau_\delta<\infty\}\right].
\end{align}
Notice that $h_{2,d_{\delta,s}}=s$ on $\{h_{d_{\delta,s}}\notin
V_\delta,d_{\delta,s}<\infty\}$, and (\ref{Plane:hitting
cone::2}) implies that we only need to prove (\ref{Plane:hitting cone::1}) when
$0<|z_1|<z_2/\delta$ and $z_2$ small enough.\medskip

\noindent \em 2. Locally, the Brownian fluctuations dominate the stable 
fluctuations. \em As $a^{-1/\alpha}S_{at}$ has the same distribution as $S_t$
for $a>0$, we have
$$\bP\left\{\theta^{1/\alpha}|S_t|\leq \frac{1}{2}\sqrt{2\kappa t\ln\ln(1/ t)}\right\}=
\bP\left\{|S_{1}|\leq \frac{1}{2}\theta^{-1/\alpha}
t^{1/2-1/\alpha} \sqrt{2\kappa\ln\ln(1/ t)}\right\}
\rightarrow 1,\ $$ when $t\downarrow 0$. Hence we can find $t_0$
such that $\bP\{\theta^{1/\alpha}|S_{t}|\leq
\frac{1}{2}\sqrt{2\kappa t\ln\ln(1/ t)}\}\geq 1/2$ for $0<t<t_0$.
Now let $s>0$ such that
 \begin{align}\label{Plane:hitting cone::3}s<t_0\wedge{2} \exp\left\{-\frac{1}{2}\exp{\frac{288}{\kappa\delta^2}}\right\}=:t_1
 \end{align}
 and
  let $z\in \mathbb{H}$ such that $0<|z_1|<s/\delta$ and $ z_2=s$.
  By (\ref{Plane:hitting cone::3}), for $0<t<s$,
  \begin{align}\label{Plane:hitting cone::4}
&\bP\left\{U_{t}\geq \sqrt{2\kappa t\ln\ln(1/ t)}/2\right\}\nonumber\\\geq&
\bP\left\{B_{t}\geq \sqrt{2t\ln\ln(1/ t)}\right\}\bP\left\{\theta^{\frac{1}{\alpha}}|S_{t}|\leq \sqrt{2\kappa t\ln\ln(1/ t)}/2\right\}\nonumber\\
\geq &\frac{1}{2}\bP\left\{B_1\geq \sqrt{2\ln\ln (1/t)}\right\}\nonumber\\
\geq&\frac{1}{4\sqrt{2\pi}\ln(1/t)\sqrt{2\ln\ln (1/t)}}.\end{align}
The last inequality of (\ref{Plane:hitting cone::4}) follows from
$\int_{x}^\infty e^{-y^2/2}\ dy\geq
 \frac{1}{2x}e^{-x^2/2}\ dy$ for $x>1$.\medskip

\noindent \em 3. $h_{2,t}$ decreases quickest if $h_{1,t}=0$, and $h_{1,t}$ 
reflects high values of $U_t$. \em By (\ref{EQ:h_1:h_2}), for each $y>0$ with 
$h_{2,0}=y$ we have 
\begin{align}\label{Plane:hitting cone::5}
h_{2,u}>y/2,\qquad\mbox{when $0<u<3y^2/16$.}
\end{align}
Therefore, if $U_{s^2/16}\geq {s}\sqrt{2\kappa\ln\ln (16/s^2)}/8$,
then by (\ref{Plane:hitting cone::3}) and
 (\ref{Plane:hitting cone::5}),
\begin{align}
|h_{1,s^2/16}|=&\bigg|z_1+\int_0^{s^2/16}\frac{2h_{1,u}}{h_{1,u}^2+h_{2,u}^2}du-U_{s^2/16}\bigg|\nonumber\\
\geq&|U_{s^2/16}|-s/\delta-\int_0^{s^2/16}\frac{2}{ s}du\nonumber\\
\geq& {s}\sqrt{2\kappa\ln\ln (4/s^2)}/8-2s/\delta\nonumber\\
\geq &s/\delta,\nonumber
\end{align}
which leads to
\begin{align}\label{Plane:hitting cone::6}
\left\{U_{s^2/16}\geq {s}\sqrt{\kappa\ln\ln (16/s^2)}/8\right\}\subseteq
\left\{\tau_\delta\leq s^2/16\right\}.
\end{align}
By (\ref{Plane:hitting cone::4}) and (\ref{Plane:hitting cone::6}),
we obtain
\begin{align}\label{Plane:hitting cone::7}
P_z\left\{\tau_\delta\leq s^2/16\right\}\geq\bP\left\{U_{s^2/16}\geq
{s}\sqrt{2\kappa\ln\ln (4/s^2)}/8\right\} \geq
\frac{1}{8\sqrt{2\pi}\ln(4/s)\sqrt{2\ln(2\ln (4/s))}}.
\end{align}
\noindent \em 4. Consider a positive starting height $s_0<t_1$ and levels 
$s_0/2^n$, $n\ge 1$. We control $\tau_\delta$ between successive levels. \em
Define $T_n=\inf\{t\geq 0:
  h_{2,t}=s_0/2^{n}\}, n\geq 1$ and $T_0=0$. Let
$p_n=P_z\{\tau_\delta\in(T_{n-1},T_n]\}$. By (\ref{Plane:hitting
cone::5}) and (\ref{Plane:hitting cone::7}) we have $$p_1\geq
\frac{1}{8\sqrt{2\pi}\ln(4/s_0)\sqrt{2\ln(2\ln (4/s_0))}}.$$ 
By the Markov property, (\ref{Plane:hitting cone::5}) and
(\ref{Plane:hitting cone::7}), we have
\begin{align}
p_n=&E_z\left[P_z\left[\tau_\delta\in(T_{n-1},T_n]\left|\mathcal{F}_{T_{n-1}}\right.\right]\right]
\nonumber\\
\geq
&E_z\left[I_{\{\tau_\delta>T_{n-1}\}}P_{h_{T_{n-1}}}\left\{|h_{1,T_{n-1}}|<s_0/(2^{n-1}\delta),\tau_\delta\leq \left(\frac{s_0}{2^{n-1}}\right)^2/16\right\}\right]
\nonumber\\
\geq & \frac{1}{8\sqrt{2\pi}(\ln(4/s_0)+(n+1)\ln2)\sqrt{2\ln(2\ln
(4/s_0)+2(n+1)\ln2)}}P_z\left\{\tau_\delta>T_{n-1}\right\}
\nonumber\\
= & \frac{1}{8\sqrt{2\pi}(\ln(4/s_0)+(n+1)\ln2)\sqrt{2\ln(2\ln
(4/s_0)+2(n+1)\ln2)}}\left(1-\sum_{k=1}^{n-1}p_k\right).\nonumber\end{align}
\noindent \em 5. We conclude. \em Now the proof is complete if we show $\sum_{n\ge 1}p_n=1$.
Otherwise, we would have $\sum_{n\ge 1}p_n<1$ and
\begin{align}
\sum_{n\ge 1}p_n\ge&
\sum_{n\ge 1}\frac{1}{8\sqrt{2\pi}(\ln(4/s_0)+(n+1)\ln2)\sqrt{2\ln(2\ln
(4/s_0)+2(n+1)\ln2)}}\left(1-\sum_{k=1}^{n-1}p_k\right) \nonumber
\\
\geq&
\sum_{n\ge 1}\frac{1}{8\sqrt{2\pi}(\ln(4/s_0)+(n+1)\ln2)\sqrt{2\ln(2\ln (4/s_0)+2(n+1)\ln2)}}\left(1-\sum_{k\ge 1}p_k\right)\nonumber\\
=&\infty,\nonumber
\end{align}which is a contradiction, so we must have $\sum_{n\ge 1}p_n=1$ as required.
\qed\medskip

\begin{Lemma}\label{Plane:from cone}
Let $z=z_1+iz_2\in \mathbb{H}$ and let $0<\kappa<4$. Then
 for any $\varepsilon>0$, there exists $\delta>0$ such that $ P_z\{\zeta<\infty\}<\varepsilon$
 for $z\in V_\delta$, the double wedge of slope $\delta$.
\end{Lemma}
\noindent{\bf Proof}$\ $
For convenience, we will  use the notation of Lemmas \ref{E:harmonic} and 2.2.
For example,  we still use notation $\tau_{a,b}$ and
$\tau_{b,a}$ for the inner and outer exit times of $(h_{1,t})_{t\geq0}$ from
$\{x\in\bR:a<|x|<b\}$. We also denote the exit time by $\tau=\tau_{a,b}\wedge\tau_{b,a}$. 
For $c\geq 0$ and a $C^2$ function $f$, set
\begin{align}\label{Plane:from cone::00**}
A_cf(y)=\frac{2y}{y^2+c^2}\partial_yf(y)+\frac{\kappa}{2} \partial^2_yf(y)+
\theta\Delta^{\alpha/2}_yf(y),\ \mbox{for $y\neq0$.}
\end{align}
Let $\beta=(2/\kappa-1/2)\wedge(1-\alpha)$ if
$\alpha<1$ and
$\beta=(2/\kappa-1/2)\wedge 1/2$ if
$1\leq\alpha<2$. Then we have $4\kappa^{-1}(1+\beta)^{-1}-1>0$. Let
$0<k<\varepsilon^{1/\beta}\wedge 1$ and let $\delta$ be a
positive number such that
\begin{align}\label{Plane:from cone::1}\delta<k\sqrt{\frac{4}{\kappa(1+\beta)}-1}
\end{align}
Define $f=w_{1-\beta}$. Noticing that
$\Delta^{\alpha/2}w_{1-\beta}(y)\leq0$, and applying
(\ref{Plane:from cone::1}),
 we have  for any $|y|>kz_1$ and  $0\leq c\leq \delta z_1$
\begin{align}\label{Plane:from cone::2}
A_cf(y)\leq&\frac{2y}{y^2+c^2}\partial_yf(y)+\frac{\kappa}{2}\partial^2_yf(y)
\nonumber\\
=&
\frac{\beta}{|y|^{2+\beta}}\left(\frac{2y^2}{y^2+c^2}-\frac{\kappa(1+\beta)}{2}\right)
\nonumber\\
\leq&
\frac{-\beta}{|y|^{2+\beta}}\left(\frac{2}{1+\delta^2/k^2}-\frac{\kappa(1+\beta)}{2}\right)
\nonumber
\\
\leq&0.\end{align} Let $\tau=\tau_{a,b}\wedge\tau_{b,a}$ for
$kz_1\leq a<z_1<b$. By Dynkin's formula,
\begin{align}
E_zf(h_{1,\tau})=&z_1^{-\beta} +E_z\int_0^{\tau}A_{h_{2,u}}f(h_{1,u-})\
du.
\end{align}
Hence by (\ref{Plane:from cone::2}) and $h_{2,u}\leq \delta z_1$, we
obtain  $E_zf(h_{1,\tau})\leq z_1^{-\beta}$. Therefore, by Remark \ref{rfu}
\begin{align}
P_z\{\zeta<\infty\} \leq&
\lim_{b\uparrow\infty}P_z\{\tau_{kz_1,b}<\tau_{b,kz_1}\}
\nonumber\\
\leq& \lim_{b\uparrow\infty}\frac{z_1^{-\beta}-\int_{\{|y|\geq b\}
}|y|^{-\beta}\ {\mu_{ b,kz_1}(dy)}}{\int_{\{|y|\leq kz_1\}}|y|^{-\beta}
{\mu_{kz_1,b}(z,dy)}-\int_{\{|y|\geq b\}}|y|^{-\beta} \ {\mu_{
b,kz_1}(z,dy)}}
\nonumber\\
\leq& k^{\beta}<\varepsilon,\nonumber
\end{align}which completes the proof for $0<\kappa<4$.
\qed

\begin{theorem}\label{thm4.3}
Let $0<\kappa < 4$. For any $z\in\overline{\mathbb{H}}\setminus
\{0\}$, we have $ P_z\{\zeta=\infty\}=1$.\end{theorem} 

\noindent{\bf
Proof}$\ $ When $z_2=0$, the conclusion follows from Proposition
\ref{Line:<4}. When $z_2>0$, the conclusion follows from Lemmas
\ref{Plane:hitting cone} and 4.2.\qed\medskip

\subsection{The supercritical phase $\kappa>4$}

We first show that we control the return time to the imaginary axis outside
an asymptotically negligible event. This will be useful when we choose 
regeneration points on the imaginary axis.

\begin{Lemma}\label{ES:Plane:>4}\begin{enumerate}\item[\rm(1)] 
Let $\kappa>4, 1\leq \alpha<2$ and let
$z=z_1+iz_2\in\mathbb{\overline{H}}\setminus \{0\}$. Denote
$\widetilde{\tau}=\inf\{t\geq0: h_{1,t-}=0\}$. Then
$\widetilde{\tau}<\infty$ with probability one. Moreover, if
$\alpha\geq 2-4/\kappa$, then there exists a constant $c$
and an event $\Theta$ such that
\begin{align}\label{ES:Plane:>4::1}
E_z\left[I_\Theta\widetilde{\tau}\right]\leq c |z_1|^{1-4/\kappa},\ \ \
P_z[\Theta^c]< c |z_1|^{1-4/\kappa},\qquad\mbox{for $0<|z_1|<1$.}
\end{align}
If $1\leq \alpha< 2-4/\kappa$, then for any $0<\beta<
1-4/\kappa $, there exists a constant $c_\beta$ and an event
$\Theta_\beta$ such that
\begin{align}\label{ES:Plane:>4::2}
 E_z\left[I_{\Theta_\beta}\widetilde{\tau}\right]
\leq c_\beta |z_1|^{\beta},\ \ \ P_z[\Theta_\beta^c]< c_\beta
|z_1|^\beta,\qquad\mbox{for $0<|z_1|<1$.}
\end{align}
Specifically we can take $\Theta$ and $\Theta_\beta$ both to be $
 \{\omega\in\Omega: \tau_{0,2}(\omega)<\tau_{2,0}(\omega)\}$ in (\ref{ES:Plane:>4::1}) and (\ref{ES:Plane:>4::2}).
\item[\rm(2)] Let $\kappa>4$ and $0<\alpha<1$, then (4.14) is true.
\end{enumerate}
\end{Lemma}
\noindent{\bf Proof}$\ $ Define $A_c$ by (\ref{Plane:from
cone::00**}). By Lemma \ref{E:harmonic}, we have $A_cw_{\beta}\leq0$
for $\beta=\alpha\wedge(2-4/\kappa)$. Then, applying the same method
as the proof of  Proposition \ref{Line:>4}, we can prove the first
conclusion.

Now let $\alpha\geq 2-4/\kappa$. By the same arguments as in
(3.5) we have
\begin{align}\label{ES:Plane:>4::3}
P_z\{\tau_{0,2}>\tau_{2,0}\}\leq
\frac{z_1^{1-4/\kappa}}{\int_{\{|y|\geq 2\} }
|y|^{1-4/\kappa}\ {\mu_{ 2,0}(z,dy)}}.
\end{align}
Let $f(x)=x^2\wedge M$ for $x\in \R$ and $M>0$. Choose $M$ big
enough such that
 $\theta\Delta^{\alpha/2}f(y)\geq -\kappa/2$ for $|y|\leq 2$. Set $\Theta=
 \{\tau_{0,2}<\tau_{2,0}\}$.
Taking  notation of Lemma \ref{Plane:from cone}, we have by
Dynkin's formula
\begin{align}\label{ES:Plane:>4::4}
E_z\left[f\left(h_{1,\tau_{0,2}\wedge\tau_{2,0}}\right)\right]\ge & z_1^2+E_z
\left[I_{\Theta}\int_0^{\widetilde{\tau}}A_{h_{2,u}}f(h_{1,u-})\ du\right]\nonumber\\
\geq&z_1^2+E_z\left[I_{\Theta}\int_0^{\widetilde{\tau}}\left(\frac{4h_{1,u-}^2}{h_{1,u-}^2+h_{2,u}^2}+\frac{\kappa}{2}\right)\ du\right]\nonumber\\
\geq&\frac{\kappa}{2} E_z\left[I_{\Theta}\widetilde{\tau}\right].
\end{align}
By (\ref{ES:Plane:>4::3}), we have
$$E_z\left[f\left(h_{1,\tau_{0,2}\wedge\tau_{2,0}}\right)\right]\leq \frac{Mz_1^{1-4/\kappa}}{\int_{\{|y|\geq 2\} }
|y|^{1-4/\kappa}\ {\mu_{ 2,0}(z,dy)}}.$$ Hence
  (\ref{ES:Plane:>4::1}) follows from (\ref{ES:Plane:>4::4}). We omit the proof of
  the other results as
  they can be proved in the same way.\qed\medskip

\begin{Lemma}\label{Lemma:sequence}
Let $\beta>0$. Let $(a_n)_{n\geq 0}$ be a sequence positive numbers
such that $a_1<(1+1/\beta)^{-1/\beta}$ and $a_{n+1}\leq
a_n-a_n^{1+\beta}/\beta$. Then
\begin{align}
a_n\leq {(a_1^{-\beta}+n-1)^{-1/\beta}},\qquad\mbox{for all $n\geq 1$.}\nonumber
\end{align}
\end{Lemma}
\noindent{\bf Proof}$\ $ It is easy to see that the assertion is
true for $n=1$. Now suppose that the assertion is true for $n=k$.
Notice that $f(x)=x+x^{\beta+1}/\beta$ is a increasing function on
$(0,(1+1/\beta)^{-1/\beta})$
 we have
$$a_{k+1}\leq a_k-a_k^{\beta+1}/\beta\leq {(a_1^{-\beta}+k-1)^{-1/\beta}}-
{(a_1^{-\beta}+k-1)^{-(\beta+1)/\beta}}/\beta\leq
{(a_1^{-\beta}+k)^{-1/\beta}},$$ which completes the
proof.\qed\medskip

\begin{theorem}\label{Th:Plane:>4}
Let $\kappa > 4$.  Then the following assertions are true.
\begin{enumerate}\item[\rm(1)] When $1\leq \alpha <2$, then for any
$z\in\overline{\mathbb{H}}\setminus \{0\}$, we have $
P_z\{\zeta<\infty\}=1$.
\item[\rm(2)] When $0<\alpha<1$, then $\lim_{|z|\downarrow 0}P_z\{\zeta<\infty
\}=1 $.
\end{enumerate}
\end{theorem}
\noindent{\bf Proof}$\ $ (1) When $z_2=0$, the conclusion follows
from Proposition \ref{Line:>4}. Next, we assume
 $z_2>0$ and, without loss of generality,  $z_1>0$. By Proposition VIII.4 in \cite{BER}, there exists a constant positive number $k_1$ such that
\begin{align}\label{Th:Plane:>4::1}
\bP\{|S_1|> x\}\leq k_1 x^{-\alpha},\qquad\mbox{for all $x>0$.}
\end{align}
  Denote $\beta=1/4-1/\kappa$.
Let $a_1$ be an arbitrary positive number such that
 \begin{align}\label{Th:Plane:>4::2}
 a_1<z_2\wedge\left(\frac{\beta}{10}\right)^{1/\beta}<\left(1+\frac{1}{\beta}\right)^{-1/\beta}.
\end{align}
  Denote $\eta_0=0$ and
$\xi_1=\inf\{t\geq 0: h_{2,t}=a_1\}$. By (\ref{EQ:h_1:h_2}), we can
check $\xi_1<\infty$ a.s.. Set
$$b_1=a_1-\frac{a_1^{1+\beta}}{\beta};\ \ \eta_1=\inf\{t\geq \xi_1:
h_{1,t}=0\}.$$ By the Markov  property and Lemma \ref{ES:Plane:>4} we
have $\eta_1<\infty$ a.s.. Define by induction
 $$a_{n+1}=h_{2,\eta_{n}};\  \xi_{n+1}=\eta_n+ \frac{5a_{n+1}^{2+\beta}}{4\beta};\
 b_{n+1}=a_{n+1}-\frac{a_{n+1}^{1+\beta}}{\beta};\ \eta_{n+1}=\inf\{t\geq \xi_{n+1}: h_{1,t}=0\}.$$
By the definitions above and  Lemma \ref{ES:Plane:>4} we see that
$\xi_n\leq \eta_n< \xi_{n+1}\leq \eta_{n+1}<\infty$, and these are sums of 
decreasing amounts of waiting time and subsequent return times of $h_t$ to the 
imaginary axis. We will show that for almost all $n\ge 1$, we have good control
of real and imaginary parts of $h_t$ so as to deduce that we reach zero in 
finite time. Specifically, set
\begin{align}\label{Th:Plane:>4::3} \ds E_n=\bigcap_{t\in[\eta_{n-1},\xi_{n}]}\{|h_{1,t}|\leq a_n\};\ \ \ \ H_n=\{
h_{2,\xi_n}\leq b_n\}.
\end{align}
Next we prove a lemma for preparation.
\begin{Lemma}\label{Th:Plane:>4:Lemma} We have
 \begin{align}\label{Th:Plane:>4::4}
\ds&P_z\left[E_n^c\left|\mathcal{F}_{\eta_{n-1}}\right.\right]\leq
\sqrt{\frac{160\kappa}{\beta\pi}}a_n^{\beta/2}\exp
\left\{{-\frac{\beta a_n^{-\beta}}{{40\kappa}}}\right\}+
\frac{10k_1\theta}{4^{1-\alpha}\beta}a_n^{{2+\beta}-\alpha};\\
&\label{Th:Plane:>4::5} E_n\subseteq H_n.
\end{align}
\end{Lemma}
\noindent{\bf Proof}$\ $ Denote $\xi_n'=\inf\{t\geq 0:h_{2,t}=
a_n/2 \}$. By (\ref{EQ:h_1:h_2}), we can prove
\begin{align}\label{Th:Plane:>4::6}h_{2,\xi_n}>a_n/2.
\end{align}
In fact, if $h_{2,\xi_n}\leq a_n/2$ we have $\xi_n'<\xi_n$ and hence
\begin{align}\label{Th:Plane:>4::7}
\frac{a_n}{2}=h_{2,\xi_n'}=&a_n+\int_{\eta_{n-1}}^{\xi_n'}\frac{-2h_{2,u}}{h_{1,u}^2+h_{2,u}^2}du\nonumber\\
\geq&a_n-\int_{\eta_{n-1}}^{\xi_n'}\frac{2}{h_{2,u}}du\nonumber\\
>&a_n-5a_n^{1+\beta}/\beta.
\end{align}
By (\ref{Th:Plane:>4::7}), we have $a_n< 10a_n^{1+\beta}/\beta\leq
10a_1^\beta a_n/\beta$, which  contradicts
(\ref{Th:Plane:>4::2}).

 By (\ref{Th:Plane:>4::2}), (\ref{Th:Plane:>4::6})
 and (\ref{EQ:h_1:h_2}), for $\eta_{n-1}
<t\leq \xi_n$, we have
\begin{align}\label{Th:Plane:>4::8}
|h_{1,t}|=&\bigg|\int_{\eta_{n-1}}^{t}\frac{2h_{1,u}}{h_{1,u}^2
+h_{2,u}^2}du+U_{t}-U_{\eta_{n-1}}\bigg|\nonumber\\
\leq&|U_{t}-U_{\eta_{n-1}}|+\int_{\eta_{n-1}}^{\xi_n}\frac{4}{a_n}\ du\nonumber\\
=&|U_{t}-U_{\eta_{n-1}}|+5a_n^{1+\beta}/\beta\nonumber\\\leq&|U_{t}-U_{\eta_{n-1}}|+a_n/2
.
\end{align}
By the reflection principle and (\ref{Th:Plane:>4::1}),
\begin{align}\label{Th:Plane:>4::9}
&P_z\left[\left.\sup_{\eta_{n-1}<t\leq
\xi_n}|U_{t}-U_{\eta_{n-1}}|>a_n/2\right|\mathcal{\eta}_{n-1}\right]
\nonumber\\ \leq&
2P_z\left[\left.\sqrt{\kappa}|B_{\xi_n}-B_{\eta_{n-1}}|>a_n/4\right|\mathcal{\eta}_{n-1}\right]
+2P_z\left[\left.\theta^{1/\alpha}|S_{\xi_n}-S_{\eta_{n-1}}|>a_n/4\right|\mathcal{\eta}_{n-1}\right]
\nonumber\\ \leq&
2P_z\left[\left.|B_1|>\beta^{1/2}a_n^{-\beta/2}/\sqrt{20\kappa}\right|\mathcal{\eta}_{n-1}\right]
+2P_z\left[\left.|S_1|>\left(\frac{4\beta}{5\theta}\right)^{1/\alpha}a_n^{1-(2+\beta)/\alpha}/4\right|\mathcal{\eta}_{n-1}\right]\nonumber\\
\leq
&\sqrt{\frac{160\kappa}{\beta\pi}}a_n^{\beta/2}\exp\left\{{-\frac{\beta a_n^{-\beta}}{{40\kappa}}}\right\}+
\frac{10k_1\theta}{4^{1-\alpha}\beta}a_n^{{2+\beta}-\alpha}.
\end{align}
Combining (\ref{Th:Plane:>4::8}) and (\ref{Th:Plane:>4::9}), we
obtain the first inequality in (\ref{Th:Plane:>4::4}).

Now suppose $ |h_{1,u}|\leq a_n$ when $\eta_{n-1}\leq u\leq \xi_{n}
$. Then we have $$\frac{2h_{2,u}}{h_{1,u}^2+a_n^2/4}>\frac{4}{5a_n}.$$
By (\ref{Th:Plane:>4::6}),
\beqs
h_{2,\xi_n}=a_n+\int_{\eta_{n-1}}^{\xi_n}\frac{-2h_{2,u}}{h_{1,u}^2+h_{2,u}^2}du\leq a_n-\int_{\eta_{n-1}}^{\xi_n}\frac{4}{5a_n}du
=a_n-a_n^{1+\beta}/\beta=b_n,
\eeqs
which proves (\ref{Th:Plane:>4::5}). \qed\medskip

\noindent\textbf{Continuation of the proof of Theorem 4.6:} Denote
\begin{align}\label{Th:Plane:>4::10}
\widetilde{\tau}_{0,n}=&\eta_n\wedge\inf\left\{t\geq \xi_n: h_{1,t}=0, \ |h_{1,u}|<2\mbox{ for }\xi_n<u<t\right\};\nonumber\\
\widetilde{\tau}_{2,n}=&\eta_n\wedge\inf\left\{t\geq \xi_n: h_{1,t}=2, \
|h_{1,u}|>0\mbox{ for }\xi_n<u<t\right\}
\end{align}
 By Lemma \ref{ES:Plane:>4}, there exists a constant $k_2>0$ such that
  \begin{align}\label{Th:Plane:>4::11}
E_z\left[\left.I_{\{\widetilde{\tau}_{0,n}<\widetilde{\tau}_{2,n}\}}\left(\eta_n-\xi_n\right)\right|\mathcal{F}_{\xi_n}\right]<
k_2 |h_{1,\xi_n}|^{1/2-2/\kappa},\ \ \
 P_z\left[\widetilde{\tau}_{0,n}>\widetilde{\tau}_{2,n}\left|\mathcal{F}_{\xi_n}\right.\right]
 < k_2 |h_{1,\xi_n}|^{1/2-2/\kappa},
\end{align}
when $0<|h_{1,\xi_n}|<1$. Denote
$F_n=\{\widetilde{\tau}_{0,n}<\widetilde{\tau}_{2,n}\}\cap E_n$ and
set $F=\bigcap_{n\geq 1}F_n$. By (\ref{Th:Plane:>4::5}) and Lemma
\ref{Lemma:sequence}
\begin{align}\label{Th:Plane:>4::12}
\bigcap_{n=1}^{N-1}E_n\subseteq\bigcap_{n=1}^{N}\left\{a_n<(a_1^{-\beta}+n-1)^{-1/\beta}\right\},\qquad\mbox{for all $N \in {\mathbb{N}}$.}
\end{align}
Write $d_n=a_1^{-\beta}+n-1$. By (\ref{Th:Plane:>4::2}), (\ref{Th:Plane:>4::3}),
(\ref{Th:Plane:>4::11}) and (\ref{Th:Plane:>4::12}),
\begin{align}\label{Th:Plane:>4::13}
P_z\left[F\right]=&\lim_{N\rightarrow\infty}P_z\left[\bigcap_{n=1}^NF_n\right]\nonumber\\
=&\lim_{N\rightarrow\infty}E_z\left[I_{\bigcap_{n=1}^{N-1}F_n}I_{E_N}P_z\left[\widetilde{\tau}_{0,N}>\widetilde{\tau}_{2,N}\left|\mathcal{F}_{\xi_N}\right.\right]\right]
\nonumber\\
\geq
&\lim_{N\rightarrow\infty}E_z\left[I_{\bigcap_{n=1}^{N-1}F_n}I_{E_N}\left(1-k_2|h_{1,\xi_N}|^{1/2-2/\kappa}\right)\right]
\nonumber\\
\geq&\lim_{N\rightarrow\infty}E_z\left[I_{\bigcap_{n=1}^{N-1}F_n}I_{E_N}
\left(1-k_2|a_N|^{1/2-2/\kappa}\right)\right]
\nonumber\\
\geq&\lim_{N\rightarrow\infty}E_z\left[I_{\bigcap_{n=1}^{N-1}F_n}I_{E_N}
\left(1-k_2d_N^{-2}\right)\right]
\nonumber\\
=&\lim_{N\rightarrow\infty}\left(1-k_2d_N^{-2}\right)E_z\left[I_{\bigcap_{n=1}^{N-1}F_n}P_z\left[
E_N\left|\mathcal{F}_{\eta_{N-1}}\right.\right]\right]
\nonumber\\
\ge&\lim_{N\rightarrow\infty}\left(1-k_2d_N^{-2}\right)E_z\left[I_{\bigcap_{n=1}^{N-1}F_n}
\left(1-\sqrt{\frac{160\kappa}{\beta\pi}}a_N^{\beta/2}\exp\left\{{-\frac{\beta a_N^{-\beta}}{{40\kappa}}}\right\}-
\frac{10k_1\theta}{4^{1-\alpha}\beta}a_N^{{2+\beta}-\alpha}\right)\right]
\nonumber\\
\ge&\lim_{N\rightarrow\infty}\left(1-k_2d_N^{-2}\right)
\left(1-\sqrt{\frac{160\kappa}{\beta\pi}}d_N^{-1/2}\exp\left\{{-\frac{\beta d_N}{{40\kappa}}}\right\}-
\frac{10k_1\theta}{4^{1-\alpha}\beta}d_N^{-1-(2-\alpha)/\beta}\right)
P_z\left[\bigcap_{n=1}^{N-1}F_n\right]
\nonumber\\
\geq&\displaystyle\prod_{n=1}^{\infty}\left(1-k_2d_n^{-2}\right)
\left(1-\sqrt{\frac{160\kappa}{\beta\pi}}d_n^{-1/2}\exp\left\{{-\frac{\beta d_n}{{40\kappa}}}\right\}-
\frac{10k_1\theta}{4^{1-\alpha}\beta}d_n^{-1-(2-\alpha)/\beta}\right)
\nonumber\\
\geq&1-\sum_{n=1}^{\infty}\left(k_2d_n^{-2}
+\sqrt{\frac{160\kappa}{\beta\pi}}d_n^{-1/2}\exp\left\{{-\frac{\beta d_n}{{40\kappa}}}\right\}+
\frac{10k_1\theta}{4^{1-\alpha}\beta}d_n^{-1-(2-\alpha)/\beta}\right).
\end{align}
By the definition of $d_n$ and (\ref{Th:Plane:>4::13}), we have
\begin{align}\label{Th:Plane:>4::14}
\lim_{a_1\downarrow 0}P_z[F]=1.
\end{align}
Set $\ds\xi=\lim_{n\rightarrow\infty}\xi_n$ and $\xi_0=0$. By
Lebesgue's monotone convergence theorem,
(\ref{Th:Plane:>4::2}), (\ref{Th:Plane:>4::11}) and
(\ref{Th:Plane:>4::12}),
\begin{align}
E_z\left[I_F\xi\right]= &\lim_{n\rightarrow \infty}E_z\left[I_F\xi_n\right]
\nonumber\\
=&
 \lim_{n\rightarrow \infty}\sum_{k=1}^nE_z\left[I_F\left(\xi_k-\eta_{k-1}\right)\right]+
 \lim_{n\rightarrow \infty}\sum_{k=1}^nE_z\left[I_F\left(\eta_{k-1}-\xi_{k-1}\right)\right]
 \nonumber\\
=&
 \lim_{n\rightarrow \infty}\sum_{k=1}^nE_z\left[I_F\frac{5a_{k}^{2+\beta}}{4\beta}\right]+
 \lim_{n\rightarrow \infty}\sum_{k=1}^nE_z\left[E_z\left[I_F(\eta_{k-1}-\xi_{k-1})\left|
 \mathcal{F}_{\xi_{k-1}}\right.\right]\right]
\nonumber\\
\leq&
 \sum_{k=1}^{\infty}E_z\left[I_F\frac{5d_{k}^{-1-2/\beta}}{4\beta}\right]+
 \sum_{k=1}^\infty E_z\left[E_z\left[\left.I_{\bigcap_{s=1}^{k-1}E_s}I_{\{\widetilde{\tau}_{0,k-1}>\widetilde{\tau}_{2,k-1}\}}\left(\eta_{k-1}-\xi_{k-1}\right)\right|
 \mathcal{F}_{\xi_{k-1}}\right]\right]
  \nonumber\\
\leq&
 \sum_{k=1}^{\infty}\frac{5d_{k}^{-1-2/\beta}}{4\beta}+
 \sum_{k=1}^\infty E_z\left[I_{\bigcap_{s=1}^{k-1}E_s}E_z\left[\left.I_{\{\widetilde{\tau}_{0,k-1}>\widetilde{\tau}_{2,k-1}\}}(\eta_{k-1}-\xi_{k-1})\right|
 \mathcal{F}_{\xi_{k-1}}\right]\right]\nonumber\\
 \leq&
 \sum_{k=1}^{\infty}\frac{5d_{k}^{-1-2/\beta}}{4\beta}+
 \sum_{k=1}^\infty E_z\left[I_{\bigcap_{s=1}^{k-1}E_s} k_2 |h_{1,\xi_{k-1}}|^{1/2-2/\kappa}\right]
  \nonumber\\
  \leq&
 \sum_{k=1}^{\infty}\frac{5d_{k}^{-1-2/\beta}}{4\beta}+
 \sum_{k=1}^\infty k_2 E_z\left[I_{\bigcap_{s=1}^{k-1}E_s}  a_{k-1}^{1/2-2/\kappa}\right]
  \nonumber\\
  \leq&
 \sum_{k=1}^{\infty}\frac{5d_{k}^{-1-2/\beta}}{4\beta}+
 \sum_{k=1}^\infty k_2 d_{k-1}^{-2}
  \nonumber\\
  <&\infty.
 \end{align}
By (\ref{Th:Plane:>4::12}), we see that $
F\subseteq\{\lim_{n\rightarrow \infty}a_n=0\}$. Hence by the definition
of $\xi$, we see $h_{2,\xi}=0$ on $F$. From this fact and
Proposition 3.2, we know $\zeta<\infty$ on $F$. Notice $a_1$ can be
 arbitrary small, we obtain the conclusion by
(\ref{Th:Plane:>4::14}).

By the same proof as  above we see that (2) can also be proved.
\qed\medskip

\subsection{Remaining critical and boundary values $\kappa=4$ and $\kappa=0$}

For $z=z_1+iz_2$ with $z_2\geq 0$, denote
\begin{align}
\widetilde{w}_p(z)=(z_1^2+z_2^2)^{(p-1)/2},\ \ \ p\neq 1;\ \ \
\widetilde{w}_1=\ln(z_1^2+z_2^2).
\end{align}
For function $f$ on the upper half plane, we set
\begin{align}\label{Plane:from cone::0}
Af(z)=\frac{-2z_2}{z_1^2+z_2^2}\partial_{z_2}f(z)+\frac{2z_1}{z_1^2+z_2^2}\partial_{z_1}f(z)+\frac{\kappa}{2}
\partial^2_{z_1}f(z)+ \theta\Delta^{\alpha/2}_{z_1}f(z).
\end{align}

\begin{Lemma}\label{Harmonic on plane}
For $0<p<\alpha+1$ and $\theta=0$,
\begin{align}\label{Harmonic on plane::0}
A\widetilde{w}_p&=\frac{p-1}{2}(z_1^2+z_2^2)^{(p-5)/2}((\kappa-4)z_2^2
+(4+\kappa(p-2))z_1^2),\nonumber\\
A\widetilde{w}_1&=(\kappa-4)(z_1^2+z_2^2)^{-2}(z_2^2-z_1^2).
\end{align}
\end{Lemma}
\noindent{\bf Proof}$\ $ When $p\neq1$, we have
\begin{align}
Af(z)=&-2(p-1)(z_1^2+z_2^2)^{{(p-5)/2}}z_2^2+2(p-1)(z_1^2+z_2^2)^{{(p-5)/2}}z_1^2\nonumber\\
&+\frac{1}{2}\kappa(p-1)(z_1^2+z_2^2)^{{(p-3)/2}}+\frac{1}{2}\kappa(p-1)(p-3)(z_1^2+z_2^2)^{{(p-5)/2}}z_1^2
\nonumber\\
=&(p-1)(z_1^2+z_2^2)^{{(p-5)/2}}(-2z_2^2+2z_1^2+\frac{\kappa}{2}(z_1^2+z_2^2)+\frac{\kappa}{2}
(p-3)z_1^2)\nonumber\\
=&\frac{p-1}{2}(z_1^2+z_2^2)^{{(p-5)/2}}((\kappa-4)z_2^2
+(4+\kappa(p-2))z_1^2).\nonumber\end{align} The second equality can
also be verified directly. \qed\medskip
\begin{Remark}\rm 
By (\ref{Harmonic on plane::0}), when $\theta=0$ we have
\begin{align}
A\widetilde{w}_{2-4/\kappa}=\frac{(\kappa-4)^2}{2\kappa}(z_1^2+z_2^2)^{-3/2-2/\kappa}z_2^2,
\end{align}
and hence $A\widetilde{w}_1=0$ for $\kappa=4$.
\end{Remark}

\begin{Lemma}\label{ES:frac on plane}For each $0<p<\alpha+1$, there exists a constant $c$ such that
\begin{align}\label{ES:frac on plane::0}
|\Delta_{z_1}^{\alpha/2}\widetilde{w}_p(z)|&\leq
c(|z_1|^{p-1-\alpha}\wedge|z_2|^{p-1-\alpha}),\ \ \ \ for \ z\neq0,
|z|<1, z\in\mathbb{\overline{H}}.
\end{align}
\end{Lemma}
\noindent{\bf Proof}$\ $ First we see the case $p<1$. We claim that
function
\begin{align}
\varphi(t):=\lim_{\varepsilon\downarrow0} \int_{\{y:|y|>\varepsilon\}}
\frac{((y+1)^2+t^2)^{(p-1)/2}-(1+t^2)^{(p-1)/2}}{|y|^{1+\alpha}}\
dy\nonumber
\end{align}
is bounded for $t\in[-1,1]$. In fact, we have for $|t|\leq1$
\begin{align}
|\varphi(t)|=& \left|\int_{-\infty}^\infty I_{\{|y|>1/2\}}
\frac{((y+1)^2+t^2)^{(p-1)/2}-(1+t^2)^{(p-1)/2}
}{|y|^{1+\alpha}}\ dy\right|\nonumber
\\+&
\left|\int_{-1/2}^{1/2} \frac{((y+1)^2+t^2)^{(p-1)/2}-
(1+t^2)^{(p-1)/2}-(p-1)(1+t^2)^{(p-3)/2}y}{|y|^{1+\alpha}}\
dy\right|\nonumber
\\\leq&
\int_{-\infty}^\infty I_{\{|y|>1/2\}}
\frac{|y+1|^{p-1}+1}{|y|^{1+\alpha}}\ dy\nonumber
\\
+&\int_{-1/2}^{1/2}
\frac{|p-1|((\frac{1}{2})^2+t^2)^{(p-3)/2}|y|^2+
|(p-1)(p-3)|(\frac{3}{2})^2((\frac{1}{2})^2+t^2)^{(p-5)/2}|y|^2}{|y|^{1+\alpha}}\
dy\nonumber
\\\leq&
\int_{-\infty}^\infty I_{\{|y|>1/2\}}
\frac{|y+1|^{p-1}+1}{|y|^{1+\alpha}}\ dy +\int_{-1/2}^{1/2}
\frac{|p-1|2^{3-p}+
|(p-1)(p-3)|(\frac{3}{2})^22^{p-5}}{|y|^{\alpha-1}}\ dy\nonumber\\
<&\infty,\nonumber
\end{align}
which gives the bound of $\varphi$ on $[-1,1]$. We denote this bound
by $c_1$. Hence for $|z_2/z_1|\leq 1$,
 we have
\begin{align}
\label{ES:frac on plane::2}
\left|\Delta_{z_1}^{\alpha/2}\widetilde{w}_p(z)\right|&=\left|\lim_{\varepsilon\downarrow0}
\mathcal{A}(1,-\alpha)\int_{\{y:|y-z_1|>\varepsilon\}}\!\!\!\!\frac{(y^2+z_2^2)^{(p-1)/2}-(z_1^2+z_2^2)^{(p-1)/2}}{|y-z_1|^{1+\alpha}}\
dy
\right|\nonumber\\
&=\mathcal{A}(1,-\alpha)|z_1|^{p-\alpha-1}\left|\lim_{\varepsilon\downarrow0}
\int_{\{y:|y-1|>\varepsilon\}}\!\!\!\!\!\!\!\!
\frac{(y^2+(z_2/z_1)^2)^{(p-1)/2}-(1+(z_2/z_1)^2)^{(p-1)/2}}{|y-1|^{1+\alpha}}\
dy
\right|\nonumber\\
&\leq  c_1 \mathcal{A}(1,-\alpha)|z_1|^{p-\alpha-1}.
\end{align}
On the other hand
\begin{align}
\label{ES:frac on plane::3}
&\left|\Delta_{z_1}^{\alpha/2}\widetilde{w}_p(z)\right|\nonumber\\
=&\mathcal{A}(1,-\alpha)|z_1|^{p-\alpha-1}\lim_{\varepsilon\downarrow0}
\left|\int_{\{y:|y|>\varepsilon\}}
\frac{((y+1)^2+(z_2/z_1)^2)^{(p-1)/2}-(1+(z_2/z_1)^2)^{(p-1)/2}}{|y|^{1+\alpha}}\
dy\right|
\nonumber\\
=&\mathcal{A}(1,-\alpha)|z_2|^{p-\alpha-1}\lim_{\varepsilon\downarrow0}
\left|\int_{\{y:|y|>\varepsilon\}}
\frac{((y+(z_1/z_2))^2+1)^{(p-1)/2}-(1+(z_1/z_2)^2)^{(p-1)/2}}{|y|^{1+\alpha}}\
dy\right|.\end{align} 
By similar calculations as above, we can also find
a positive number $c_2$ such that
\begin{align}
\label{ES:frac on plane::4} \lim_{\varepsilon\downarrow0}
\left|\int_{\{y:|y|>\varepsilon\}}
\frac{((y+(z_1/z_2))^2+1)^{(p-1)/2}-(1+(z_1/z_2)^2)^{(p-1)/2}}{|y|^{1+\alpha}}\
dy\right|\leq c_2
\end{align}
for $|z_1/z_2|<1$. Combining (\ref{ES:frac on
plane::2}), (\ref{ES:frac on plane::3}) and (\ref{ES:frac on plane::4}), we
get
\begin{align}
&\left|\Delta_{z_1}^{\alpha/2}\widetilde{w}_p(z)\right|\leq (c_1+c_2)
\mathcal{A}(1,-\alpha)(|z_1|^{p-\alpha-1}\wedge|z_2|^{p-\alpha-1})\nonumber\end{align}
which completes the proof for $p<1$. The case $p\geq 1$ can be
checked with the same method.\qed\medskip

\begin{theorem}\label{Th:Plane:=4}
Let $\kappa = 4$. Then for any $z\in\overline{\mathbb{H}}\setminus
\{0\}$, we have $ P_z\{\zeta=\infty\}=1$.
\end{theorem}
\noindent{\bf Proof}$\ $ As in the  case of the real line, we need
to construct a continuous  function $f$ which is subharmonic with
respect to $A$ on a pointed neighbourhood of zero and satisfies
\begin{align}\label{Plane:=4::1}\lim_{|z|\downarrow0}f(z)=-\infty;\ \ \ \lim_{|z|\uparrow\infty}f(z)\geq 0.
\end{align}
First we see the case $\alpha>1$. Let $f_1$ be a continuous function
on $\overline{\mathbb{H}}$ such that
$$f_1(z)=-\widetilde{w}_{2-\alpha/2},
\ \ |z|\leq 1,z\in\mathbb{\overline{H}};\ \ \ \ f_1(z)=0,\ \
|z|>2,z\in\mathbb{\overline{H}}.$$ By (\ref{ES:frac on plane::0}) we
can check that there exists a positive number $c_1$ such that
\begin{align}\label{Plane:=4::2}
\left|\Delta_{z_1}^{\alpha/2}f_1(z)\right|&\leq
c_1\left(|z_1|^{1-3\alpha/2}\wedge|z_2|^{1-3\alpha/2}\right),\qquad\mbox{for $|z|<1/2$, $z\in\mathbb{\overline{H}}$.}
\end{align}
By (\ref{Harmonic on plane::0}) and (\ref{ES:frac on plane::0}),
there exist positive numbers $c_2$ and $c_3$ such that
\begin{align}\label{Plane:=4::3}
Af_1(z)&\geq c_2(z_1^2+z_2^2)^{-(\alpha+2)/4},\qquad\mbox{for $\theta=0$ and 
$z\in\mathbb{\overline{H}}$.}
\end{align}
and
\begin{align}\label{Plane:=4::4}
|\Delta_{z_1}^{\alpha/2}\widetilde{w}_1(z)|&\leq
c_3(|z_1|^{-\alpha}\wedge|z_2|^{-\alpha}),\ \ \ \
z\in\mathbb{\overline{H}}.
\end{align}
Denote $f=f_1+\widetilde{w}_{1}$. It is easy to see that $f$
satisfies (\ref{Plane:=4::1}). By
(\ref{Plane:=4::2}), (\ref{Plane:=4::3}), (\ref{Plane:=4::4}), and
noticing that $-(\alpha+2)/2<-\alpha<1-3\alpha/2$, we
get
$$\lim_{|z|\downarrow0}Af(z)=\infty.$$  Hence by (2) in Lemma 2.3 and Dynkin's formula we finish the proof of  $\alpha>1$.
When $0<\alpha\leq 1$, the proof is still valid provided that we
define $f_1$ by
$$f_1(z)=-\widetilde{w}_{1+\alpha/2},\
 \ |z|\leq 1,z\in\mathbb{\overline{H}};\ \ \ \ f_1(z)=0,\ \ |z|>2,z\in\mathbb{\overline{H}}.$$
When $\theta=0$ we can simply choose $f=\widetilde{w}_1$.
\qed\medskip

Next we consider the pure jump case, i.e. $\kappa=0$. The proof for
this case is similar to   the case of $0<\kappa<4$. For
$\delta,\gamma>0$, denote $V_{\gamma,\delta}=\{z=(z_1,z_2):0<z_2\leq
\delta|z_1|^{\gamma/2}\}$ and $\sigma_{\gamma,\delta}=\inf\{t\geq 0:
h_t\in V_{\gamma,\delta}\}$.

\begin{Lemma}\label{Plane:hitting cone,kappa=0}
If $\kappa=0$ and $0<\alpha<2$, then for each $\delta>0$ and
$z\in{\mathbb{H}}$,
\begin{align}\label{Plane:hitting cone,kappa=0::1}
P_z\{\sigma_{\alpha,\delta}<\infty\}=1.
\end{align}\end{Lemma}
\noindent{\bf Proof}$\ $ We only need to prove
 (\ref{Plane:hitting cone,kappa=0::1}) when $z\notin V_{\alpha,\delta}$. Without loss of
generality we assume that $\delta<1$. By arguments similar to the
case of $0<\kappa<4$, we only need to prove (\ref{Plane:hitting
cone,kappa=0::1}) when $0<|z_1|^{\alpha/2}<z_2/\delta$ and
$z_2$ small enough.

Now let $s>0$ such that
 \begin{align}\label{Plane:hitting cone::3''}
 s<4 \exp\left\{-\frac{1}{2}\exp\left\{3\left(2^{4/\alpha}\right)\delta^{-2/\alpha}\theta^{-1/\alpha}\right\}\right\}=:t_1
 \end{align}
 and
  let $z\in \mathbb{H}$ such that $0<|z_1|^{\alpha/2}<s/\delta$ and $ z_2=s$.
   By Proposition VIII.4 in \cite{BER}, there exists a positive number $k_1$ such that for $0<t<s$,
  \begin{align}\label{Plane:hitting cone::4''}
&\bP\left\{U_{t}\geq {\left(\theta t\right)^{1/\alpha}\ln\ln(1/ t)}\right\}=
\bP\left\{S_{1}\geq \ln\ln(1/ t)\right\}\geq k_1\left(\ln\ln(1/t)\right)^{-\alpha}.
\end{align}
We claim that if $U_{s^2/16}\geq2^{-4/\alpha}
\theta^{1/\alpha}s^{2/\alpha} {\ln\ln (16/s^2)}$,
then
\begin{align}\label{Plane:hitting cone::4'}
|h_{1,u}|\geq (s/\delta)^{2/\alpha},\qquad\mbox{for some $u\in(0,s^2/16]$.}
\end{align}
If this is not true,
 by (\ref{Plane:hitting cone::5}) and (\ref{Plane:hitting
 cone::3''}),
\begin{align}
\left|h_{1,s^2/16}\right|=&\left|z_1+\int_0^{s^2/16}\frac{2h_{1,u}}{h_{1,u}^2+h_{2,u}^2}du-U_{s^2/16}\right|\nonumber\\
\geq&\left|U_{s^2/16}\right|-(s/\delta)^{2/\alpha}-\int_0^{s^2/16}\frac{8(s/\delta)^{2/\alpha}}{ s^2}du\nonumber\\
\geq& 2^{-4/\alpha}
\theta^{1/\alpha}s^{2/\alpha} {\ln\ln (16/s^2)}-2(s/\delta)^{2/\alpha}\nonumber\\
\geq & (s/\delta)^{2/\alpha},\nonumber
\end{align}
which leads to a contradiction. By (\ref{Plane:hitting cone::4'})
\begin{align}\label{Plane:hitting cone::6''}
\left\{U_{s^2/16}\geq 2^{-4/\alpha}
\theta^{1/\alpha}s^{2/\alpha} {\ln\ln
(16/s^2)}\right\}\subseteq \left\{\sigma_{\alpha,\delta}\leq s^2/16\right\}.
\end{align}
By (\ref{Plane:hitting cone::4''}) and (\ref{Plane:hitting
cone::6''}), we obtain
\begin{align}\label{Plane:hitting cone::7''}
P_z\left\{\sigma_{\alpha,\delta}\leq s^2/16\right\}\geq\bP\left\{U_{s^2/16}\geq
2^{-4/\alpha} \theta^{1/\alpha}s^{2/\alpha}
{\ln\ln (16/s^2)}\right\} \geq k_1({\ln\ln (16/s^2)})^{-\alpha}\end{align}
Let $s_0$ be a  positive number such that  $s_0<t_1/4$. Define 
$T_n=\inf\{t\geq 0:
  h_{2,t}=s_0/2^{n}\}, n\geq 1$ and $T_0=0$. Let
$p_n=P_z\{\sigma_{\alpha,\delta}\in(T_{n-1},T_n]\}$.  By the Markov property,
(\ref{Plane:hitting cone::5}) and (\ref{Plane:hitting cone::7''}),
we have
\begin{align}
p_n=&E_z\left[P_z\left[\sigma_{\alpha,\delta}\in(T_{n-1},T_n]\left|\mathcal{F}_{T_{n-1}}\right.\right]\right]
\nonumber\\
\geq
&E_z\left[I_{\{\sigma_{\alpha,\delta}>T_{n-1}\}}P_{h_{T_{n-1}}}\left\{|h_{1,T_{n-1}}|^{\alpha/2}<s_0/(2^{n-1}\delta),\sigma_{\alpha,\delta}\leq \left(\frac{s_0}{2^{n-1}}\right)^2/16\right\}\right]
\nonumber\\
\geq &k_1(\ln(2(n+1)\ln2-2\ln
s_0)^{-\alpha}P_z\{\sigma_{\alpha,\delta}>T_{n-1}\}
\nonumber\\
\geq & k_1(\ln(2(n+1)\ln2-2\ln
s_0)^{-\alpha}\left(1-\sum_{k=1}^{n-1}p_k\right),\nonumber\end{align} Hence we
can prove (\ref{Plane:hitting cone,kappa=0::1}) by the same method
as in the case of $0<\kappa<4$.\qed\medskip

Recall that we denote $\tau_{a,b}=\inf\{t>0: h_{1,t}\leq a; h_{1,u}<b,\mbox{ for all
$0\leq u<t$}\}$.

\begin{Lemma}\label{Plane:from cone''}
Let $z=(z_1,z_2)\in \mathbb{H}\setminus \{0\}$,  $\kappa=0$.
\begin{enumerate}\item[\rm(1)] If  $0<\alpha\leq 1$, then $ P_z\{\zeta<\infty\}=0$.
\item[\rm(2)] If $1< \alpha<2$, for any $\varepsilon>0$, there
exists $\delta>0$ such that
  $ P_z\{0,\tau_{0,c(\theta,\alpha)}<\tau_{c(\theta,\alpha),0}\}<\varepsilon$
 for $z$ satisfying $0<|z_2|/|z_1|^{\alpha/2}<\delta$ and $0<|z_1|<
 c(\theta,\alpha):=(2
 \mathcal{A}(1,-\alpha)\gamma
(\alpha,\frac{1}{2})\theta)^{-1/(2-\alpha)}$.
\end{enumerate}
\end{Lemma}
\noindent{\bf Proof}$\ $
 For convenience, we will  use the notation of Lemma \ref{Plane:from cone}. Here we set
\begin{align}\label{Plane2:from cone::0}
A_cf(y)=\frac{2y}{y^2+c^2}\partial_yf(y)+
\theta\Delta^{\alpha/2}_yf(y),\qquad\mbox{for $y\in\bR\setminus\{0\}$.}
\end{align}
for any $C^2$ function $f$. When $0<\alpha<1$,  we can check that
$A_cw_{(\alpha+1)/2}(y)<0$ for $y\neq 0$. We can also check
that $A_cw_1(y)\geq 0$ for $y\neq 0$. Hence we can prove (1) by
Dynkin's formula.

Next we assume $1< \alpha<2$.  Let $0<|z_1|<c(\theta,\alpha)$. For
any $\varepsilon>0$, let $0<k<\varepsilon^{2}\wedge 1$ and let
$\delta$ be a positive number such that
\begin{align}\label{Plane:from cone::1''}\delta<\left(\frac{k^\alpha}{2
\mathcal{A}(1,-\alpha)\gamma (\alpha,\frac{1}{2})\theta
}\right)^{1/2}
\end{align}
Define $f=w_{1/2}$. We claim that   $A_c f<0$ if
\begin{align}
\label{Plane:from cone::1.5''} k|z_1|<|y|<c(\theta,\alpha),\ \ \
0\leq c\leq \delta |z_1|^{\alpha/2}.
 \end{align} In fact when
$k^2|z_1|^2<|y|^2<\delta^2|z_1|^\alpha$, by (\ref{Plane:from
cone::1''})
\begin{align}\label{Plane:from cone::2''}
A_cf(y)=&\frac{-
|y|^{1/2}}{y^2+c^2}+\mathcal{A}(1,-\alpha)\gamma\left(\alpha,\frac{1}{2}\right)\theta
|y|^{-1/2-\alpha}
\nonumber\\
\leq & |y|^{-1/2-\alpha}\left(\frac{- |y|^\alpha}{y^2+\delta^2
|z_1|^\alpha}+\mathcal{A}(1,-\alpha)\gamma\left(\alpha,\frac{1}{2}\right)\theta\right)
\nonumber\\
\leq& |y|^{-1/2-\alpha}\left(\frac{- k^\alpha}{2\delta^2
}+\mathcal{A}(1,-\alpha)\gamma\left(\alpha,\frac{1}{2}\right)\theta\right)\nonumber
\\
\leq&0.\end{align} Similarly,  when
$c(\theta,\alpha)^2>|y|^2\geq\delta^2|z_1|^\alpha$,
\begin{align}\label{Plane:from cone::3''}
A_cf(y)\leq & |y|^{-1/2-\alpha}\left(\frac{-
|y|^\alpha}{2y^2}+\mathcal{A}(1,-\alpha)\gamma\left(\alpha,\frac{1}{2}\right)\theta\right)\leq
0.
\end{align}
Combing (\ref{Plane:from cone::2''}) and (\ref{Plane:from
cone::3''}), we get the claim. Thus, applying  Dynkin's formula to
$f$, we have
\begin{align}
P_z\left\{\tau_{0,c(\theta,\alpha)}<\infty\right\} \leq&
P_z\left\{\tau_{k|z_1|,c(\theta,\alpha)}<\tau_{c(\theta,\alpha),k|z_1|}\right\}
\nonumber\\
\leq&\frac{|z_1|^{-1/2}-\int_{\{|y|\geq c(\theta,\alpha)\}
}|y|^{-1/2}\ {\mu_{
c(\theta,\alpha),k|z_1|}(z,dy)}}{\int_{\{|y|\leq k|z_1|\}
}|y|^{-1/2}
{\mu_{k|z_1|,c(\theta,\alpha)}(z,dy)}-\int_{\{|y|\geq c(\theta,\alpha)\}
}|y|^{-1/2} \ {\mu_{ c(\theta,\alpha),k|z_1|}(z,dy)}}
\nonumber\\
\leq& k^{1/2}<\varepsilon,\nonumber
\end{align}which completes the proof.
\qed

\begin{theorem}\label{thm4.13}
Let $\kappa=0$ and $0<\alpha<2$. For any
$z\in\overline{\mathbb{H}}\setminus \{0\}$, we have $
P_z\{\zeta=\infty\}=1$.\end{theorem} 

\noindent{\bf Proof}$\ $ When
$z_2=0$, the conclusion follows from Lemma \ref{hitting, lemma }.
When $z_2>0$ and $0<\alpha\leq 1$, the conclusion follows from Lemma
\ref{Plane:hitting cone,kappa=0} and \ref{Plane:from cone''}.

Next we assume  $1<\alpha<2$ and $z\in \mathbb{H}$. For any $n\in
\mathbb{N}$ and $\varepsilon>0$, by Lemma \ref{Plane:from cone''},
there exists $\delta_n>0$ such that
$P_z\{\tau_{0,c(\theta,\alpha)}<\tau_{0,c(\theta,\alpha)}\}<\varepsilon/2^n$
for $0<|z_1|<c(\theta,\alpha)$. For any $z\in \mathbb{H}$, define
$\tau_1=\inf\{t>0; h_t\in V_{\delta_n,\alpha}\}$ and
$\sigma_1=\inf\{t\geq\tau_1; |h_{1,t}|>c(\theta,\alpha)\}$. Define
by induction, $\tau_n=\inf\{t\geq \sigma_{n-1}; h_t\in
V_{\delta_n,\alpha}, |h_{1,t}|<c(\theta, \alpha)/2\}$ and
$\sigma_n=\inf\{t\geq\tau_{n}; |h_{1,t}|>c(\theta,\alpha)\mbox{ or }
h_{t-}=0\}$ for $n\geq 2$. By Lemmas \ref{Plane:hitting
cone,kappa=0} and \ref{Plane:from cone''} as well as the quasi-left continuity
of paths, we have
$$
P_z\{\zeta<\infty\}=\sum_{n=1}^\infty P_z\{\sigma_n=\zeta<\infty\}+
P_z\left[\bigcap_{n=1}^\infty \{\sigma_n<\zeta<\infty\}\right]\leq \sum_{n=1}^\infty
\frac{\varepsilon}{2^n}=\varepsilon,
$$
which completes the proof. \qed\medskip

\subsection{Proofs of Theorem \ref{thm1} and Corollary \ref{cor1}}

The statement of Theorem \ref{thm1} is contained in Theorems \ref{thm4.3},
\ref{Th:Plane:>4}, \ref{Th:Plane:=4} and \ref{thm4.13}. To prove Corollary
\ref{cor1}, we just note that the generator of the stable process with all
jumps of size exceeding $c$ removed has as its generator
\beqs \Delta_{x|c}^{\alpha/2}w(x)=\lim_{\varepsilon\downarrow 0}\mathcal{A}(1,-\alpha)\int_{\{y:\varepsilon<|y-x|<c\}}\frac{w(y)-w(x)}{|x-y|^{1+\alpha}}dy,
\eeqs
and a computation as in Lemma \ref{E:harmonic} shows that
\beqs \Delta_{x|c}^{\alpha/2}w_p(x)=\mathcal{A}(1,-\alpha)|x|^{p-1-\alpha}\left(\gamma(\alpha,p)-\frac{p-1}{\alpha}\int_{1-x/c}^{1+x/c}v^{p-2}|1-v|^{\alpha-p}dv\right)
\eeqs
and for $x$ small enough, the right-most factor has the same sign as 
$\gamma(\alpha,p)$. It can now be checked that all arguments can be adapted.\qed

\section{The increasing cluster of SLE driven by $U=\sqrt{\kappa}B+\theta^{-1/\alpha}S$}

Denote the life time of $(h_t(z))_{t\ge 0}$ starting at $h_0(z)=z\in
\overline{\mathbb{H}}$ by $\zeta(z)$ as in Section 2.2 and define
$$K_t=\{z\in \overline{\mathbb{H}},\ \ \zeta(z)\leq t\}, $$
the associated family of strictly
increasing compact sets in $\mathbb{H}$, and $\mathbb{H}\setminus
K_t$ the associated simply connected open set. First note that unlike the Brownian case, 
$K_t$ is not always connected by the following lemma.

\begin{proposition}
$$\bP\{\mbox{$K_t$ is a disconnected set in $\overline{\bH}$}\}>0,\qquad\mbox{for all $t>0$.}$$
\end{proposition}
\noindent{\bf Proof}$\ $ Let $t>0$. Set $\tau=\inf\{s\geq 0:\
|U_s|>1\}$. By (2.4) we have for $u<\tau$
\begin{align}
|h_u(z)|=|z+\int_{0}^u\frac{2}{h_s-U_s}ds|\geq
|z|-\int_{0}^u\frac{2}{|h_s|-1}ds.\nonumber
\end{align}
Hence we can check that
\begin{align}
K_{\tau-}\subseteq B(0,2t+2),\qquad\mbox{for $\tau<t$.}
\end{align}
Denote Loewner's conformal mapping associated with $K_{\tau}$ by $g_{\tau}$, 
and $$B=\left\{U_\tau-U_{\tau-}>2\sup\left\{\left|g_{1,\tau}\left(z\right)\right|:z\in B\left(0,2t+2\right)\right\}+\left(4t+5\right)\right\}.$$ 
By (5.1), we have
\begin{align}
B\subseteq\{\mbox{$K_\tau$ is a disconnected set}\}.\end{align} Set
$B'=\{|U_s-U_\tau|\leq1,\ \tau<s<\tau+t\}$. By similar arguments 
as for (5.1) we have
\begin{align}
g_{\tau}(B(0,2t+2))\cap B(U_\tau,2t+2)=\varnothing\qquad\Rightarrow\qquad \overline{K}_{\tau-}\cap \overline{K_t\setminus K_{\tau-}}=\varnothing.
%\{z:z_1=2t+5/2\}\cap (K_t(\omega)\setminus K_{\tau-}(\omega))=\varnothing,\qquad
%\mbox{for all $\omega\in B\cap B'$.}
\end{align}
As $\bP[B\cap B']=\bP[B]\bP[B']>0$, by (5.1)-(5.3), we get the
conclusion.  \qed\medskip

\noindent{\bf Proof of Theorem \ref{thm1b}}$\ $
In what follows we denote Lebesgue measure on $\overline{\bH}$
by $m(\cdot)$. Recall that Theorem \ref{thm1b} claims the following:
(1)\ When $\kappa\leq 4$,  we have $m(\bigcup_{t>0} K_t)=0$, a.s..\ \
(2)\ When
$\kappa> 4$ and $1\leq\alpha<2$, we have $m(\mathbb{\overline{H}}\setminus
\bigcup_{t>0} K_t)=0$, a.s.. (3)\ When $\kappa> 4$ and $0<\alpha<1$, we have
$\lim_{r\downarrow 0}m(B(0,r) \cap (\bigcup_{t>0} K_t))/m(B(0,r))=1$,
a.s. and $\lim_{r\uparrow \infty}m(B(0,r) \cap (\bigcup_{t>0}
K_t))/m(B(0,r))=0$ a.s..

First we show that the lifetime function $\zeta(\omega,z)$ is 
measurable from
$(\Omega\times\mathbb{\overline{H}},\mathcal{F}\otimes
\mathcal{B}(\mathbb{\overline{H}}))$ to $(
[0,\infty],\mathcal{B}([0,\infty]))$. Denote
$\tau^z_a=\inf\{t\geq 0: h_t(z)\in B(0,a)\}$ for $h_0(z)=z$ and $a>0$. For
any $r>0$, we have
$$\{(\omega,z):\zeta(\omega,z)\le r\}=\bigcup_{k=1}^\infty\bigcap_{l=1}^\infty\{(\omega,z): z\in
\overline{\mathbb{H}}, |z|>1/k, \tau^z_{1/l}(\omega)\le r\}.$$ Hence we only
need to show that $\{(\omega,z): z\in \overline{\mathbb{H}}, |z|>a,
\tau^z_b(\omega)\le r\}\in\mathcal{F}\otimes \mathcal{B}(\mathbb{\overline{H}})$
for any $a>b>0$. As the coefficient function of the stochastic
differential equation $(2.4)$ is Lipschitz and satisfies the linear
growth condition outside any neighbourhood of zero, by Theorem
6.4.3 in \cite{APP}, we know that $(h_{t}(z))_{t\ge 0}$, $z\in
\overline{\mathbb{H}}$, have the flow property before hitting
$B(0,b)$. Therefore we have $\{(\omega,z): z\in \overline{\mathbb{H}},
|z|>a, \tau^z_b(\omega)<r\}\in\mathcal{F}\otimes
\mathcal{B}(\mathbb{\overline{H}})$.

 Now let $\kappa\leq 4$. By
Theorem \ref{thm1}(i), we have
\begin{align}
&\bE[m(\{z:
\zeta(z)<\infty\})]=\bE\left[\int_{\overline{\mathbb{H}}}I_{\{\zeta(z)<\infty\}}m(dz)\right]
\nonumber\\=&\int_{\overline{\mathbb{H}}}\bE[I_{\{\zeta(z)<\infty\}}]m(dz)=
\int_{\overline{\mathbb{H}}}P_z\{\zeta<\infty\}m(dz)=0,
\end{align}
which leads to (1). Similarly, by Theorem \ref{thm1}(ii), when $\kappa>4$ and
$1\leq\alpha<2$, we have for any $n>0$
\begin{align}
&\bE[m(\{z: \zeta(z)<\infty\},
|z|<n)]=\bE\left[\int_{\overline{\mathbb{H}}}I_{\{|z|<n\}}I_{\{\zeta(z)<\infty\}}m(dz)\right]
\nonumber\\=&\int_{\overline{\mathbb{H}}}I_{\{|z|<n\}}\bE[I_{\{\zeta(z)<\infty\}}]m(dz)=m(\{z:|z|<n\}).\nonumber
\end{align}
Hence, we have $m(\mathbb{\overline{H}}\setminus \bigcup_{t>0} K_t)=0$,
a.s.. (3) can be proved by Theorem \ref{thm1}(iii) and the same method.
\qed

 \section{$\beta$-SLE driven by $\alpha$-stable processes}

Let $(S_t)_{t\geq0}$ be the standard symmetric  $\alpha-$stable L\'{e}vy
process. For simplicity we take $(S_t)_{t\geq 0}$ as the standard
Brownian motion when $\alpha=2$. For $1<\beta\leq 2$ define
  the following generalized SLE
$(g_t)_{t\geq0}$, which we call $\beta$-SLE:
\begin{align}
\partial_tg(z)=&\frac{2|g_t(z)-\theta^{1/\alpha}S_t|^{2-\beta}}{g_t(z)-\theta^{1/\alpha}S_t},\qquad g_0(z)=z,\quad z\in\mathbb{\overline{H}}\setminus\{0\},\qquad 1<\beta\leq2,\ 0<\alpha\le 2;\nonumber
%\\{g}'_{0}(t,z)=&\frac{2|{g}_{\beta}(t,z)-\theta^{\frac{1}{\alpha}}S_t|\ln|{g}_{\beta}(t,z)-\theta^{\frac{1}{\alpha}}S_t|}
%{{g}_{\beta}(t,z)-\theta^{\frac{1}{\alpha}}S_t},\ \ \
%{g}_{0}(0,z)=z,\ z\in \mathbb{\overline{H}}\setminus\{0\}.\nonumber
\end{align}
where the derivative above is the right derivative as $S_t$ is right 
continuous. Let $h_t(z)=g_t(z)-\theta^{1/\alpha}S_t$, then
we have
\begin{align} \label{equation}
 dh_t(z)=&\frac{2|h_t(z)|^{2-\beta}}{ h_t(z)}dt -\theta^{1/\alpha}dS_t,\qquad
h_0(z)=z,\quad z\in \mathbb{\overline{H}}\setminus\{0\}.%\\
%dh_{0}(t,z)=&\frac{2|h_{0}(t,z)|\ln|h_{0}(t,z)|}{ h_{0}(t,z)}dt
%+\theta^{\frac{1}{\alpha}}dS_t,\ \ h_{0}(0,z)=z,\ z\in
%\mathbb{\overline{H}}\setminus\{0\}.
\end{align}
Here $(h_t(z))_{t\ge 0}$ is again a well defined stochastic process
up to hitting zero. In fact, similar
to the SLE model we could use a much more general driving process in the above stochastic differential  
equation. In our setting, when $x\in\R$, $(h_t(x))_{t\geq0}$ is
an $\R$-valued Markov process  and its  generator ${A}^{\alpha,\beta,\theta}$ acting
on $C^2$ function $f$ is
\begin{align}
A^{\alpha,\beta,\theta}
f(y)=&\frac{|y|^{2-\beta}}{y}\partial_yf(y)+\theta\Delta^{\alpha/2}_yf(y),\qquad
\mbox{for all $y\neq 0$, $1<\beta\leq 2$.}\label{genalphasle2}%\\ \ \ \ A^0
%f(y)=&sign(y)\ln |y|D_yf(y)+\theta\Delta^{\alpha/2}f(y),\qquad\mbox{for all $y\neq 0$,}
\end{align}% where $sign(y)=y/|y|$ for $y\neq 0$.
We also denote simply $h_t=h_t(x)$, where $h_0=x$ under $P_x$. Also the lifetime of $h_t$ is again denoted by $\zeta$.

\begin{proposition} \label{Line: GSLE}
Let $\theta>0$, $1<\beta<2$, and $x\in \R$ with $ x\neq 0$. The
following statements are valid:
\begin{enumerate}\item[\rm (a)] If $\alpha>\beta$, then 
$\limsup_{|x|\downarrow0}P_x\{\zeta=\infty\}|x|^{-\delta}<\infty$ and
$\limsup_{|x|\uparrow \infty}P_x\{\zeta<\infty\}|x|^{\delta}<\infty$ for all
$0<\delta<\alpha-1$.
\item[\rm(b)] If $\alpha=\beta$, there is a phase transition at $\theta_0(\alpha)=2/(\mathcal{A}(1,-\alpha)|\gamma(\alpha,1)|)$ as follows
  \beqs  P_x(\zeta<\infty)=1\quad\mbox{if}\quad\theta> \theta_0(\alpha)\qquad\mbox{and}\qquad
  P_x(\zeta=\infty)=1\quad\mbox{if}\quad 0<\theta\leq\theta_0(\alpha).
  \eeqs
\item[\rm(c)] If $\alpha<\beta$, then $P_x(\zeta=\infty)=1$.
\end{enumerate}
\end{proposition}
\noindent{\bf Proof}$\ $(a)\ Let $0<\delta<\alpha-1$. By Lemma
\ref{E:harmonic} we can find a positive constant $c_1$ such that
$A^{\alpha,\beta,\theta}w_{1+\delta}(y)<0$ if $0<|y|<c_1$. Hence for $0<a<x<c_1$
we have
\begin{align}P_x\{\zeta=\infty\}&\leq \lim_{a\downarrow 0}P_x\{\tau_{a,c_1}>\tau_{c_1,a}\}\nonumber\\
&\leq \lim_{a\downarrow 0}\frac{\int_{\{|y|\leq a\}}|y|^{\delta}\
{\mu_{ c_1,a}(x,dy)}-x^{\delta}}{\int_{\{|y|\leq a\}} |y|^{\delta} {\mu_{
a,c_1}(x,dy)}-\int_{\{|y|\geq c_1\}} |y|^{\delta}\ {\mu_{
c_1,a}(x,dy)}}\nonumber\\
&=x^{\delta}\left/\lim_{a\downarrow 0}\int_{\{|y|\geq c_1\}} |y|^{\delta}\
{\mu_{ c_1,a}(x,dy)}\right.,
\end{align}
which gives the first conclusion in (a). Again by Lemma
\ref{E:harmonic} we can find a positive constant $c_2$ such that
$A^{\alpha,\beta,\theta}w_{1-\delta}(y)<0$ if $|y|>c_2$. Similarly we have for
$0<c_2<x<b$
\begin{align}
P_x\{\zeta<\infty\}&\leq \lim_{b\uparrow
\infty}P_x\{\tau_{b,c_2}>\tau_{c_2,b}\} \leq
x^{-\delta}\left/\lim_{b\uparrow \infty}\int_{|y|\leq c_2 }\right.
|y|^{-\delta}\ {\mu_{ c_2,b}(x,dy)},
\end{align}
which gives the second conclusion in (a).
\medskip

(b) Let $\beta=\alpha$. Define the function
$$\varphi(p)=\frac{2(1-p)}{\mathcal{A}(1,-\alpha)\gamma(\alpha,p)},
\quad p\neq 1\qquad\mbox{and}\qquad
\varphi(1)=\frac{2}{\mathcal{A}(1,-\alpha)|\gamma(\alpha,1)|}=\theta_0(\alpha).$$ 
By Lemma \ref{E:harmonic}, we can check that $\varphi$ is a strictly
increasing continuous function on $(0,\alpha)$ and
\begin{align}\label{GSLE::Line::1} \varphi(0+):=\lim_{p\downarrow 0
}\varphi(p)>0;\qquad\lim_{p\uparrow \alpha}\varphi(p)=\infty.
\end{align} Denote by $\varphi^{-1}$ the inverse function of $\varphi$ on
$(\varphi(0+),\infty)$. By Lemma \ref{E:harmonic} and (\ref{genalphasle2}) we have
$A^{\alpha,\beta,\theta}w_{\varphi^{-1}(\theta)}=0$ for
$\theta\in(\varphi(0+),\infty)$. Hence when
$\theta\in(\varphi(0+),\infty)$, with the help of harmonic function
$w_{\varphi^{-1}(\theta)}$ we can prove the conclusion by the same
method as in Section 3. When $\theta\in (0, \varphi(0+)]$ we can
check that $A^{\alpha,\beta,\theta}w_1>0$, which also leads to our conclusion.
\medskip

(c) By Lemma 2.1 we can find a positive constant $c_3$ such that
$A^{\alpha,\beta,\theta}w_0-c_3w_0<0$.  We can prove (c) by this fact and the
same method as in Case 2 of Proposition \ref{Line:<4}.\qed\medskip

The behaviour in (a) is new. It did not occur in the same way for SLE since 
Brownian forcing is at the same time at the top of the self-similarity range 
$\alpha\in(0,2]$ and the critical forcing where the phase transition occurs,
in particular, where in the upper phase the force is strong enough to 
overcome the potential of the singularity of $h_t$ at zero. 
For $\beta$-SLE driven by an $\alpha$-stable process with $\alpha>\beta$, the 
forcing is more than just strong enough to overcome the singularity at zero, 
but on the other hand, the outward drift is stronger and makes $h_t$ transient,
so that there is positive probability that $h_t$ does not hit zero. In this,
there are similarities with $\kappa>4$ and transient driving force for SLE.

If $\alpha=2>\beta$, this can only happen if $\bR\cap\bigcup_{t\ge 0}K_t=[a,b]$
for some $-\infty<a<0<b<\infty$. This means that the $\beta$-SLE cluster then 
grows more in the vertical direction, whereas adding a transient driving force
to SLE yields clusters that grow more in the horizontal direction (and 
necessarily by disconnecting jumps).  

In what follows we concentrate on the critical and as such most interesting 
case $\beta=\alpha$. We will
show that  the phase transition indicated in Proposition \ref{Line: GSLE} 
can be extended from $z=x\in\bR$ to $z\in\bH$ in strong analogy to the 
well-known $\kappa=4$ phase transition. 
 Recall for
$\delta>0$, we denote by $V_\delta=\{z=z_1+iz_2:0<z_2\leq
\delta|z_1|)\}$ the double wedge of slope $\delta$ and by 
$\tau_\delta=\inf\{t\geq 0: h_t\in V_\delta\}$ the first entrance time of $h$.
\begin{Lemma}\label{GSLE:: <}
Let $\theta>0$. Then for each $\delta>0$ and $z\in{\mathbb{H}}$,
\begin{align}\label{GSLE:: <::1}
P_z\{\tau_{\delta}<\infty\}=1.
\end{align}\end{Lemma}

\noindent{\bf Proof}$\ $ By arguments similar to the case of Lemma
4.1, we only need to prove (\ref{Plane:hitting cone,kappa=0::1})
when $0<|z_1|<z_2/\delta$ and $z_2$ small enough.
By (\ref{equation}), for each $y>0$ with $h_{2,0}=y$ we have
\begin{align}\label{Plane2:hitting cone::5}
h_{2,u}>y/2,\qquad\mbox{when}\quad 0<u<y^\alpha/2^{2+\alpha}.
\end{align} Now let $s>0$ such that
 \begin{align}\label{Plane:hitting cone::3''**}
 s<16^{1/\alpha} \exp\left\{-\frac{1}{2}\exp\left\{3\cdot 2^{4/\alpha}\delta\theta^{-1/\alpha}\right\}\right\}=:t_1
 \end{align}
 and
  let $z\in \mathbb{H}$ such that $0<|z_1|<s/\delta$ and $ z_2=s$.

We claim that if $S_{s^\alpha/16}\geq2^{-4/\alpha} s
 {\ln\ln (16/s^\alpha)}$,
then
\begin{align}\label{Plane:hitting cone::4'**}
|h_{1,u}|\geq s/\delta,\qquad\mbox{for some $u\in (0,s^\alpha/16]$.}
\end{align}
If this is not true,
 by (\ref{Plane2:hitting cone::5}) and (\ref{Plane:hitting
 cone::3''**}),
\begin{align}
|h_{1,s^\alpha/16}|=&\left|z_1+\int_0^{s^\alpha/16}\frac{2h_{1,u}}{(h_{1,u}^2+h_{2,u}^2)^{\alpha/2}}du
-\theta^{1/\alpha}S_{s^\alpha/16}\right|\nonumber\\
\geq&\left|\theta^{1/\alpha} S_{s^\alpha/16}\right|-s/\delta-\int_0^{s^\alpha/16}\frac{2^{1+\alpha}
}{ s^{\alpha-1}\delta}du\nonumber\\
\geq& 2^{-4/\alpha}
\theta^{1/\alpha}s{\ln\ln (16/s^\alpha)}-2s/\delta\nonumber\\
\geq & s/\delta,\nonumber
\end{align}
which leads to a contradiction. By (\ref{Plane:hitting cone::4'**})
\begin{align}\label{Plane:hitting cone::6''**}
\left\{S_{s^\alpha/16}\geq 2^{-4/\alpha} s {\ln\ln
(16/s^\alpha)}\right\}\subseteq \left\{\tau_\delta\leq s^\alpha/16\right\}.
\end{align}
By (\ref{Plane:hitting cone::4''}) and (\ref{Plane:hitting
cone::6''**}), we obtain
\begin{align}\label{Plane:hitting cone::7''**}
P_z\left\{\tau_\delta\leq s^\alpha/16\right\}\geq P\left\{U_{s^\alpha/16}\geq
2^{-4/\alpha} \theta^{1/\alpha}s {\ln\ln
(16/s^\alpha)}\right\} \geq k_1\left({\ln\ln
(16/s^\alpha)}\right)^{-\alpha}.\end{align} 
Let $s_0$ be a  positive
number such that  $s_0< t_1$. Define
$T_n=\inf\{t\geq 0:
  h_{2,t}=s_0/2^{n}\}, n\geq 1$ and $T_0=0$. Let
$p_n=P_z\{\tau_\delta\in(T_{n-1},T_n]\}$.  By the Markov property,
(\ref{Plane2:hitting cone::5}) and (\ref{Plane:hitting cone::7''**}),
we have
\begin{align}
p_n=&E_z\left[P_z\left[\tau_\delta\in(T_{n-1},T_n]\left|\mathcal{F}_{T_{n-1}}\right.\right]\right]
\nonumber\\
\geq
&E_z\left[I_{\{\tau_\delta>T_{n-1}\}}P_{h_{T_{n-1}}}\left\{\left|h_{1,T_{n-1}}\right|^{\alpha/2}<s_0/\left(2^{n-1}\delta\right),
\tau_\delta\leq \left(\frac{s_0}{2^{n-1}}\right)^\alpha/16\right\}\right]
\nonumber\\
\geq &k_1(\ln(\alpha(n-1)\ln2+4\ln 2-\alpha\ln
s_0)^{-\alpha}P_z\{\tau_\delta>T_{n-1}\}
\nonumber\\
\geq & k_1(\ln(\alpha(n-1)\ln2+4\ln 2-\alpha\ln
s_0)^{-\alpha}\left(1-\sum_{k=1}^{n-1}p_k\right),\nonumber\end{align} 
Hence we
can complete the proof by the same arguments as in Lemma
4.1.\qed\medskip

\begin{proposition}\label{GSLE:: <::}
Let $1<\alpha<2$ and $0<\theta <
\theta_0(\alpha)$. For any
$z\in\overline{\mathbb{H}}\setminus \{0\}$, we have $
P_z\{\zeta=\infty\}=1$.\end{proposition} \noindent{\bf Proof}$\ $
When $z_2=0$, the conclusion follows from Proposition \ref{Line:
GSLE}. When $z_2>0$, by Lemma \ref{GSLE:: <} we only need to prove
that,
 for any $\varepsilon>0$, there exists $\delta>0$ such that $ P_z\{\zeta<\infty\}<\varepsilon$
 for $z$ satisfying $0<|z_2|/|z_1|<\delta$.
  For $c\geq 0$ and  $C^2$
function $f$, set
\begin{align}\label{Plane2:from cone::00}
A_c^{\alpha,\theta}f(y)=\frac{2y}{(y^2+c^2)^{\alpha/2}}\partial_yf(y)+
\theta\Delta^{\alpha/2}_yf(y),\quad\mbox{for $y\neq0$.}
\end{align}
Let  $\theta\in (0,\theta_0(\alpha))$ and define
$$b= \varphi^{-1}\left(\frac{\theta_0(\alpha)+(\theta\vee \varphi(0+))}{2}
\right).$$ By the
definition of $\varphi$, we see that
   $0<b<1$. Set $\theta_1=\theta/\varphi(b)$. It is easy to see
that $\theta_1<1$.
   Let
$0<k<\varepsilon^{1/(1-b)}\wedge 1$ and let $\delta$ be a
positive number such that
\begin{align}\label{GSLE:from
cone::1}\delta<k\sqrt{\theta_1^{-2/\alpha}-1}.
\end{align}
Define $f=w_{b}$  and applying (\ref{GSLE:from cone::1}),
 we have  for any $|y|>k|z_1|$ and  $0\leq c\leq \delta |z_1|$
\begin{align}\label{GSLE:from cone::2}
A_c^{\alpha,\theta}f(y)=&\frac{2(b-1)|y|^{b-1}}{(y^2+c^2)^{\alpha/2}}+\theta\mathcal{A}(1,-\alpha)
\gamma(\alpha,b) |y|^{b-1-\alpha}
\nonumber\\
\leq&
\frac{b-1}{|y|^{\alpha+1-b}}\left(\frac{2}{(1+\delta^2/k^2)^{\alpha/2}}+\theta\mathcal{A}(1,-\alpha)
\gamma(\alpha,b)/(b-1)\right)
\nonumber\\
=&
\frac{b-1}{|y|^{\alpha+1-b}}\left(\frac{2}{(1+\delta^2/k^2)^{\alpha/2}}-2\theta
/\varphi(b)\right)\nonumber
\\
=&\frac{b-1}{|y|^{\alpha+1-b}}\left(\frac{2}{(1+\delta^2/k^2)^{\alpha/2}}-
2\theta_1\right)\nonumber\\
\leq& 0 .\end{align} By (\ref{GSLE:from cone::2}) and  the same
calculation as in Lemma \ref{Plane:from cone} we have
\begin{align}
P_z\{\zeta<\infty\} \leq k^{1-b}<\varepsilon,\nonumber
\end{align}which completes the proof.
\qed\medskip

Next we consider the case $\theta>\theta_0(\alpha)$. First we prepare
a  result   corresponding to Lemma \ref{ES:Plane:>4}.
\begin{Lemma}\label{GSLE:Plane:>4} Let $1<\alpha<2$ and $\theta >
\theta_0(\alpha)$. Let
$z=(z_1,z_2)\in\mathbb{\overline{H}}\setminus \{0\}$. Denote
$\widetilde{\tau}=\inf\{t\geq0: h_{1,t-}=0\}$. Then
$\widetilde{\tau}<\infty$ with probability one. Moreover, there
exist a constant $c$ and an event $\Theta$ such that
\begin{align}\label{GSLE:Plane:>4::1}
E_z[I_\Theta\widetilde{\tau}]< c |z_1|^{\varphi^{-1}(\theta)-1},\ \
\ P_z[\Theta^c]< c |z_1|^{\varphi^{-1}(\theta)-1},\qquad\mbox{for $0<|z_1|<1$.}
\end{align}
Specifically we can take $\Theta$  to be $
 \{\tau_{0,2}<\tau_{2,0}\}$ in (\ref{GSLE:Plane:>4::1}).
\end{Lemma}
\noindent{\bf Proof}$\ $ We omit the proof as it is the same as
for Lemma \ref{ES:Plane:>4}.\qed\medskip

\begin{Lemma}\label{GSLE:Plane:>4,*} Let $1<\alpha<2$ and $\theta >
\theta_0(\alpha)$.  Let $\delta>0$ be
such that
${(\varphi^{-1}(\theta)-1)(1-\delta/\alpha)-2\delta}=:r>0$. Then there
exists   a constant number $k_3$, depending on $\alpha$, $\delta$
and $\theta$, such that for any $a>0$
\begin{align}\label{GSLE:Plane:>4::1*}
P_z\{L< a^{\alpha+\delta}/\delta\}\leq k_3a^{2\delta},\qquad\mbox{where}
\quad
L=\int_0^{3a^r}I_{\{|h_{1,t}|<a\}}\
dt.\end{align}
\end{Lemma}
\noindent{\bf Proof}$\ $
It is obvious that we can also assume  $a$ to be small enough such
that
\begin{align}\label{GSLE:Plane:>4::1**}
32a^\delta <\delta,\ \ \
a^{{(\varphi^{-1}(\theta)-1)(1-\delta/\alpha)}-2\delta}>a^{\alpha+\delta}/\delta.
\end{align}
Denote $\tau(s)=\inf\{t:t\geq s, h_{1,t}=0\}-s$ for $s>0$. By
(\ref{GSLE:Plane:>4::1}), we have
\begin{align}\label{GSLE:Plane:>4::3*}
&P_z\{|h_{1, a^{\alpha+ \delta}/\delta}|<a^{1-\delta/\alpha},
\tau({{a^{\alpha+ \delta}/\delta}})\geq
a^{{(\varphi^{-1}(\theta)-1)(1-\delta/\alpha)}-2\delta}\}\nonumber\\
\leq & c a^{2\delta}+ca^{(\varphi^{-1}(\theta)-1)(1-\delta/\alpha)}\nonumber\\
\leq &2c a^{2\delta}.
\end{align}
We claim that
\begin{align}\label{Th:Plane:>4::8***}
&\left\{\sup_{0<t\leq a^{\alpha+ \delta}/\delta}|h_{1,t}|\geq
a\right\}\subseteq \left\{\sup_{0<t\leq a^{\alpha+
\delta}/\delta}\theta^{1/\alpha}|S_{t}|\geq a/8\right\}
\\
\label{Th:Plane:>4::8****}&\left\{\sup_{0<t\leq a^{\alpha+
\delta}/\delta}|h_{1,t}|\geq a^{1-\delta/\alpha}\right\}\subseteq
\left\{\sup_{0<t\leq a^{\alpha+
\delta}/\delta}\theta^{1/\alpha}|S_{t}|\geq
a^{1-\delta/\alpha}/8\right\}
\end{align}
Let $t'=\inf\{t: |h_{1,t}(w)|\geq a\}$, $t''=\sup\{t\leq t':
|h_{1,t}(w)|< a/2\}$ and suppose that $\omega$ belongs to the left hand
side of (\ref{Th:Plane:>4::8***}), then by the first inequality of
(\ref{GSLE:Plane:>4::1**})
\begin{align}\label{Th:Plane:>4::8**}
a/2\leq
|h_{1,t'}-h_{1,t''-}|=&\left|\int_{t''}^{t'}\frac{2h_{1,u}}{(h_{1,u}^2
+h_{2,u}^2)^{\alpha/2}}du-\theta^{1/\alpha}S_{t'}+\theta^{1/\alpha}S_{t''-}\right|\nonumber\\
\leq&\left|\theta^{1/\alpha} (S_{t'}-S_{t''-})\right|+\int_{t''}^{t'}2h_{1,u}^{1-\alpha}\ du\nonumber\\
\leq&\left|\theta^{1/\alpha}(S_{t'}-S_{t''-})\right|+8a^{1+\delta}/\delta\nonumber\\\leq&\left|\theta^{1/\alpha}
(S_{t'}-S_{t''-})\right|+a/4 .
\end{align}
which proves (\ref{Th:Plane:>4::8***}). We omit the proof of
(\ref{Th:Plane:>4::8****}) as the proof is  the same. By the reflection
principle we have
\begin{align}\label{GSLE:Plane:>4::8*} \bP\left\{\sup_{0<t\leq
a^{\alpha+ \delta}/\delta}\theta^{1/\alpha}|S_{t}|\geq a/8\right\}
\leq& 2\bP\left\{|S_{a^{\alpha+ \delta}/\delta}|\geq
\theta^{-1/\alpha}a/8\right\}
\nonumber\\ \leq&2\bP\left\{|S_1|\geq \delta^{1/\alpha}a^{-\delta/\alpha}/8 \right\}\nonumber\\
\leq &2^{1+3\alpha}k_1\theta\delta^{-1} a^{ \delta}.
\end{align}
Similarly we have
\begin{align}\label{GSLE:Plane:>4::9*} &\bP\left\{\sup_{0<t\leq
a^{\alpha+ \delta}/\delta}\theta^{1/\alpha}|S_{t}|\geq
a^{1-\delta/\alpha}/8\right\} \leq 2^{1+3\alpha}k_1
\theta\delta^{-1}a^{2\delta}.
\end{align}
By (\ref{GSLE:Plane:>4::1**}), (\ref{GSLE:Plane:>4::3*}),
(\ref{Th:Plane:>4::8***}), (\ref{Th:Plane:>4::8****}),
 (\ref{GSLE:Plane:>4::8*}) and (\ref{GSLE:Plane:>4::9*}),
\begin{align}\label{Th:Plane:>4::9*}
&P_z\left\{L< a^{\alpha+\delta}/\delta\right\}\nonumber\\\leq& P_z\left\{a\leq
\sup_{0<t\leq a^{\alpha+
\delta}/\delta}|h_{1,t}|<a^{1-\delta/\alpha},L<
a^{\alpha+\delta}/\delta\right\}+P_z\left\{\sup_{0<t\leq a^{\alpha+
\delta}/\delta}|h_{1,t}|\geq a^{1-\delta/\alpha}\right\} \nonumber\\
\leq& P_z\left\{a\leq \sup_{0<t\leq a^{\alpha+
\delta}/\delta}|h_{1,t}|<a^{1-\delta/\alpha}, \tau({{a^{\alpha+
\delta}/\delta}})<a^{{(\varphi^{-1}(\theta)-1)(1-\delta/\alpha)}-2\delta},L<
a^{\alpha+\delta}/\delta\right\}\nonumber\\
&+P_z\left\{\sup_{0<t\leq a^{\alpha+
\delta}/\delta}|h_{1,t}|<a^{1-\delta/\alpha}, \tau({{a^{\alpha+
\delta}/\delta}})\geq
a^{{(\varphi^{-1}(\theta)-1)(1-\delta/\alpha)}-2\delta}\right\}+2^{1+3\alpha}k_1\theta\delta^{-1}
a^{2\delta}
\nonumber\\
\leq& P_z\left\{a\leq \sup_{0<t\leq a^{\alpha+
\delta}/\delta}|h_{1,t}|<a^{1-\delta/\alpha}, \tau({{a^{\alpha+
\delta}/\delta}})<a^{{(\varphi^{-1}(\theta)-1)(1-\delta/\alpha)}-2\delta},\right.
\nonumber\\ 
&\qquad\left.\sup_{\tau({{a^{\alpha+ \delta}/\delta}})\leq t\leq
\tau({{a^{\alpha+ \delta}/\delta}})+a^{\alpha+ \delta}/\delta}|h_{1,t}|\geq a\right\} + 2^{1+3\alpha}\delta^{-1}(k_1\theta+c)a^{2\delta}\nonumber\\
\leq& P_z\left\{\sup_{0<t\leq a^{\alpha+ \delta}/\delta}|h_{1,t}|\geq a,
\sup_{\tau({{a^{\alpha+ \delta}/\delta}})\leq t\leq
\tau({{a^{\alpha+ \delta}/\delta}})+a^{\alpha+ \delta}/\delta}|h_{1,t}|\geq a\right\} + 2^{1+3\alpha}\delta^{-1}(k_1\theta+c)a^{2\delta}\nonumber\\
\leq&
2^{1+3\alpha}\delta^{-1}(k_1\theta+c+2\delta^{-1}c^2)a^{2\delta},
\end{align}
which completes the proof. \qed\medskip

\begin{proposition}\label{GSLE:Plane:>4**}
 Let $1<\alpha<2$ and $\theta >\theta_0(\alpha)$. Let
$z\in\mathbb{\overline{H}}\setminus \{0\}$. Then  $
P_z\{\zeta<\infty\}=1$.
\end{proposition}
\noindent{\bf Proof}$\ $ The proof will follow the arguments for
Theorem \ref{Th:Plane:>4} with some technical differences.  Fix $z=z_1+iz_2\in\overline{\bH}$. When
$z_2=0$, the conclusion follows from Proposition \ref{Line: GSLE}.
Next, we assume
 $z_2>0$ and, without loss of generality,  $z_1>0$.
  Denote $\beta>0$ small enough such that
  \begin{align}\label{GSLE:Plane:>4::2*}
(\varphi^{-1}(\theta)-1)(1-\beta/\alpha)&\geq 6\beta,\\
\label{GSLE:Plane:>4::2**}\frac{{(\varphi^{-1}(\theta)-1)(1-\beta/\alpha)}-2\beta}{2\alpha}(\varphi^{-1}(\theta)-1)&\geq
2\beta.
\end{align}
Write $\widetilde{\alpha}=
{{(\varphi^{-1}(\theta)-1)(1-\beta/\alpha)}-2\beta}$. Let $a_1$ be
an arbitrary positive number such that
 \begin{align}\label{GSLE:Plane:>4::2}
 a_1<z_2\wedge\left(\frac{\beta}{\beta+1}\right)^{1/\beta}\qquad\mbox{and}\qquad a_1^{1+\beta}/\beta<a_1/2.
\end{align}
  Denote $\eta_0=0$ and
$\xi_1=\inf\{t\geq 0: h_{2,t}=a_1\}$.  Set
$$b_1=a_1-\frac{a_1^{1+\beta}}{\beta};\ \ \eta_1=\inf\{t\geq \xi_1:
h_{1,t}=0\}.$$ By Lemma \ref{GSLE:Plane:>4} we have $\eta_1<\infty$
a.s..  Define by induction
 $$a_{n+1}=h_{2,\eta_{n}};\  \xi_{n+1}=\eta_n+ {3a_{n+1}^{\widetilde{\alpha}}};\
 b_{n+1}=a_{n+1}-\frac{a_{n+1}^{1+\beta}}{\beta};\ \eta_{n+1}=\inf\{t\geq \xi_{n+1}: h_{1,t}=0\}.$$
 Let $L_n=\int_{\eta_{n-1}}^{\xi_n}I_{\{|h_{1,t}|<a_n\}}\
dt$. Define events
\begin{align}\label{Th:Plane2:>4::3} \ds E_n=\left\{L_n\geq 2^{\alpha/2}a_n^{\alpha +\beta}/\beta\right\};\qquad
G_n=\left\{|h_{1,\xi_n}|>8a_n^{\widetilde{\alpha}/2\alpha}\right\};\qquad H_n=\left\{h_{2,\xi_n}\leq b_n\right\}.
\end{align}
Next we prove the following assertions:
 \begin{align}\label{Th:Plane:>4::4*}
&G_n\subseteq \left\{
\theta^{1/\alpha}\sup_{\eta_{n-1}<t<\xi_n}
|S_{\xi_n}-S_{t}|> a_n^{\widetilde{\alpha}/2\alpha}\right\},\\
&\label{Th:Plane:>4::5*} E_n\subseteq H_n,\\
\label{Th:Plane:>4::6*}\ds&P_z[E_n^c\cup
G_n|\mathcal{F}_{\eta_{n-1}}]\leq \left(6\theta
k_1+2^{\alpha\beta/(\alpha+\beta)}k_3\right)a_n^{2\beta}.
\end{align}

Suppose that $\theta^{1/\alpha}\sup_{\eta_{n-1}<t<\xi_n}
|S_{\xi_n}-S_{t}|\leq a_n^{\widetilde{\alpha}/2\alpha}$, we
will check (\ref{Th:Plane:>4::4*}) by proving that
$|h_{1,\xi_n}|\leq 8 a_n^{\widetilde{\alpha}/2\alpha}$.
Otherwise we can find $t'\in (\eta_{n-1}, \xi_n)$ such that
$|h_{1,t'-}|\leq a_n^{\widetilde{\alpha}/2\alpha}$  and
$|h_{1,t}|\geq a_n^{\widetilde{\alpha}/2\alpha}$ for $t\in
(t',\xi_n)$. So we have
\begin{align}\label{Th:Plane2:>4::8}
|h_{1,\xi_n}|=&\bigg|\int_{\eta_{n-1}}^{\xi_n}
\frac{2h_{1,u}}{(h_{1,u}^2
+h_{2,u}^2)^{\alpha/2}}du-\theta^{1/\alpha}
S_{\xi_n}+\theta^{1/\alpha}S_{\eta_{n-1}}\bigg|\nonumber\\
\leq &\bigg|\int_{t'}^{\xi_n} \frac{2h_{1,u}}{(h_{1,u}^2
+h_{2,u}^2)^{\alpha/2}}du-\theta^{1/\alpha}
S_{\xi_n}+\theta^{1/\alpha}S_{t'-}\bigg|+ |h_{1,t'-}|\nonumber\\
\leq&\bigg|\int_{t'}^{\xi_n}2h_{1,u}^{1-\alpha}du\bigg|+2
a_n^{\widetilde{\alpha}/2\alpha}\nonumber\\
\leq&6a_n^{\widetilde{\alpha}(1+{\alpha})/2\alpha}+2
a_n^{\widetilde{\alpha}/2\alpha}\nonumber\\
\leq&8a_n^{\widetilde{\alpha}/2\alpha}.
\end{align}

Now suppose that $L_n\geq 2^{\alpha/2}a_n^{\alpha
+\beta}/\beta$. If $h_{2,\xi_n}<a_n/2$, by the second inequality of
(\ref{GSLE:Plane:>4::2}), we see that (\ref{Th:Plane:>4::5*}) is
true. When $h_{2,\xi_n}\geq a_n/2$, we have
\begin{align}
h_{2,\xi_n}=&a_n+\int_{\eta_{n-1}}^{\xi_n}\frac{-2h_{2,u}}{(h_{1,u}^2+h_{2,u}^2)^{\alpha/2}}du\nonumber\\
\leq&a_n-\int_{\eta_{n-1}}^{\xi_n}\frac{a_n}{(h_{1,u}^2+a_n^2)^{\alpha/2}}du\nonumber\\
\leq&a_n-2^{-\alpha/2}\int_{\eta_{n-1}}^{\xi_n}I_{\{|h_{1,t}|<a_n\}}a_n^{1-\alpha}du\nonumber\\
\leq&a_n-a_n^{1+\beta}/\beta=b_n,\nonumber
\end{align}
which completes the proof of (\ref{Th:Plane:>4::5*}).
 (\ref{Th:Plane:>4::6*}) can be proved by Lemma
 \ref{GSLE:Plane:>4,*}, (\ref{Th:Plane:>4::4*}), (\ref{Th:Plane:>4::5*})
 and the following results.
\begin{align}\label{GSLE:Plane:>4::9}
&P_z\left[\left.\theta^{1/\alpha}\sup_{\eta_{n-1}<t<\xi_n}
|S_{\xi_n}-S_{t}|>a_n^{\widetilde{\alpha}/2\alpha}
\right|\mathcal{F}_{\eta_{n-1}}\right]\leq
2P_z\left[\left.|S_{\xi_n-\eta_{n-1}}|>\theta^{-1/\alpha}a_n^{
\widetilde{\alpha}/2\alpha}\right|\mathcal{F}_{\eta_{n-1}}\right]\nonumber\\
\leq& 2P_z\left[\left.|S_1|>3^{-1/\alpha}\theta^{-1/\alpha}
a_n^{-\widetilde{\alpha}/2\alpha}\right|\mathcal{F}_{\eta_{n-1}}\right]\leq
6\theta k_1 a_n^{\widetilde{\alpha}/2}\leq 6\theta k_1 a_n^{2\beta},
\end{align}
where we used (\ref{GSLE:Plane:>4::2*}) in the last inequality of
(\ref{GSLE:Plane:>4::9}).

 As for SLE we denote
\begin{align}\label{Th:Plane2:>4::10}
\widetilde{\tau}_{0,n}=&\inf\{t\geq \xi_n: h_{1,t}=0, \ |h_{1,u}|<2\mbox{ for }\xi_n<u<t\};\nonumber\\
\widetilde{\tau}_{2,n}=&\inf\{t\geq \xi_n: h_{1,t}\geq 2, \
|h_{1,u}|>0\mbox{ for }\xi_n<u<t\}
\end{align}
 By Lemma \ref{GSLE:Plane:>4}, there exists a constant $k_4>0$ such that
  \begin{align}\label{Th:Plane2:>4::11}
E_z\left[\left.I_{\{\widetilde{\tau}_{0,n}<\widetilde{\tau}_{2,n}\}}(\eta_n-\xi_n)\right|\mathcal{F}_{\xi_n}\right]<
k_4 |h_{1,\xi_n}|^{\varphi^{-1}(\theta)-1},\ \ \
 E_z\left[\left.I_{\{{\widetilde{\tau}_{0,n}>\widetilde{\tau}_{2,n}}\}}\right|\mathcal{F}_{\xi_n}\right]
 < k_4 |h_{1,\xi_n}|^{\varphi^{-1}(\theta)-1},
\end{align}
when $0<|h_{1,\xi_n}|<1$. Denote
$F_n=\{\widetilde{\tau}_{0,n}<\widetilde{\tau}_{2,n}\}\cap (E_n\cap
G_n^c)$ and set $F=\bigcap_{n\geq 1}F_n$. By (\ref{Th:Plane:>4::5*})
and Lemma \ref{Lemma:sequence}
\begin{align}\label{Th:Plane2:>4::12}
\bigcap_{n=1}^{N-1}(E_n\cap
G_n^c)\subseteq\bigcap_{n=1}^{N}\left\{a_n<\left(a_1^{-\beta}+n-1\right)^{-1/\beta}\right\},\
\ \ \forall N \in {\mathbb{N}}.
\end{align}Write $d_n=a_1^{-\beta}+n-1$. By (\ref{GSLE:Plane:>4::2**}), (\ref{Th:Plane:>4::6*}),
(\ref{Th:Plane2:>4::11}) and (\ref{Th:Plane2:>4::12}),
\begin{align}\label{Th:Plane2:>4::13}
P_z[F]=&\lim_{N\rightarrow\infty}P_z\left[\bigcap_{n=1}^NF_n\right]\nonumber\\
=&\lim_{N\rightarrow\infty}E_z\left[I_{\bigcap_{n=1}^{N-1}F_n}I_{E_N\cap
G_N^c}P_z\left[\widetilde{\tau}_{0,N}<\widetilde{\tau}_{2,N}\left|\mathcal{F}_{\xi_N}\right.\right]\right]
\nonumber\\
\geq
&\lim_{N\rightarrow\infty}E_z\left[I_{\bigcap_{n=1}^{N-1}F_n}I_{E_N\cap
G_N^c}\left(1-k_4 \left|h_{1,\xi_N}\right|^{\varphi^{-1}(\theta)-1}\right)\right]
\nonumber\\
\geq&\lim_{N\rightarrow\infty}E_z\left[I_{\bigcap_{n=1}^{N-1}F_n}I_{E_N\cap
G_N^c}
\left(1-2^{3(\varphi^{-1}(\theta)-1)}k_4\left|a_N\right|^{(\varphi^{-1}(\theta)-1)\widetilde{\alpha}/2\alpha}\right)\right]
\nonumber\\
\geq&\lim_{N\rightarrow\infty}E_z\left[I_{\bigcap_{n=1}^{N-1}F_n}I_{E_N\cap
G_N^c} \left(1-2^{3(\varphi^{-1}(\theta)-1)}k_4d_N^{-2}\right)\right]
\nonumber\\
=&\lim_{N\rightarrow\infty}\left(1-2^{3(\varphi^{-1}(\theta)-1)}k_4d_N^{-2}\right)E_z\left[I_{\bigcap_{n=1}^{N-1}F_n}P_z\left[
E_N\cap G_N^c\left|\mathcal{F}_{\eta_{N-1}}\right.\right]\right]
\nonumber\\
\geq&\lim_{N\rightarrow\infty}\left(1-2^{3(\varphi^{-1}(\theta)-1)}k_4d_N^{-2}\right)
\left(1-\left(6\theta k_1+2^{2\beta/(\alpha+\beta)}k_3\right)a_N^{2\beta}\right)
P_z\left[\bigcap_{n=1}^{N-1}F_n\right]
\nonumber\\
\geq&\displaystyle\prod_{n=1}^{\infty}\left(1-2^{3(\varphi^{-1}(\theta)-1)}k_4d_n^{-2}\right)
\left(1-\left(6\theta k_1+2^{2\beta/(\alpha+\beta)}k_3\right)d_n^{-2}\right)
\nonumber\\
\geq&1-\sum_{n=1}^{\infty}\left(6\theta
k_1+2^{2\beta/(\alpha+\beta)}k_3+2^{3(\varphi^{-1}(\theta)-1)}k_4\right)d_n^{-2}
.
\end{align}
By the definition of $d_n$ and (\ref{Th:Plane2:>4::13}), we have
\begin{align}\label{Th:Plane2:>4::14}
\lim_{a_1\downarrow 0}P_z[F]=1.
\end{align}
 By
Lebesgue's monotone convergence  theorem,
(\ref{GSLE:Plane:>4::2**}), (\ref{Th:Plane2:>4::11}) and
(\ref{Th:Plane2:>4::12}),
\begin{align}
E_z\left[I_F\zeta\right]= &\lim_{n\rightarrow \infty}E_z\left[I_F\xi_n\right]
\nonumber\\
=&
 \lim_{n\rightarrow \infty}\sum_{k=1}^nE_z\left[I_F(\xi_k-\eta_{k-1})\right]+
 \lim_{n\rightarrow \infty}\sum_{k=1}^nE_z\left[I_F(\eta_{k-1}-\xi_{k-1})\right]
 \nonumber\\
\leq&
 \sum_{k=1}^{\infty}3E_z\left[I_F{d_k^{\widetilde{\alpha}/\beta}}\right]+
 \sum_{k=1}^\infty E_z\left[E_z\left[\left.I_{\bigcap_{s=1}^{k-1}(E_s\cap G_s^c)}
 I_{\{{\widetilde{\tau}_{0,k-1}>\widetilde{\tau}_{2,k-1}}\}}(\eta_{k-1}-\xi_{k-1})\right|
 \mathcal{F}_{\xi_{k-1}}\right]\right]
 \nonumber\\
 \leq&
 \sum_{k=1}^{\infty}3{d_k^{\widetilde{\alpha}/\beta}}+
 \sum_{k=1}^\infty E_z\left[I_{\bigcap_{s=1}^{k-1}(E_s\cap G_s^c)} 2^{3(\varphi^{-1}(\theta)-1)}k_4
 \left|h_{1,\xi_{k-1}}\right|^{(\varphi^{-1}(\theta)-1)\widetilde{\alpha}/2\alpha}\right]
  \nonumber\\
  \leq&
 \sum_{k=1}^{\infty}3{d_k^{\widetilde{\alpha}/\beta}}+
 \sum_{k=1}^\infty 2^{3(\varphi^{-1}(\theta)-1)}k_4E_z\left[I_{\bigcap_{s=1}^{k-1}(E_s\cap G_s^c)}  a_{k-1}^{(\varphi^{-1}(\theta)-1)\widetilde{\alpha}/2\alpha}\right]
  \nonumber\\
  \leq&
 \sum_{k=1}^{\infty}3{d_{k}^{-4}}+
 \sum_{k=1}^\infty2^{3(\varphi^{-1}(\theta)-1)}k_4d_{k-1}^{-2}
  \nonumber\\
  <&\infty,\nonumber
 \end{align}
which completes the proof.\qed\medskip

%We denote $K_t$ to be the increasing cluster of $\alpha$-SLE.\medskip

\noindent{\bf Proofs of Theorem \ref{thm2} and Corollary \ref{cor2}}$\ $
The statement of Theorem \ref{thm2} is contained in Propositions 
\ref{GSLE:: <::} and \ref{GSLE:Plane:>4**}.
The proof of the corollary is the same as for SLE with the help
of these propositions.\qed
\medskip

\subsection*{Acknowledgements}

The second author would like to thank Wendelin Werner for introducing him to
SLE, Terry Lyons for asking what happens if you replace the driving Brownian
motion by a L\'evy process. Both authors would like to thank them for 
discussions and continued interest in this work. 

\bibliographystyle{abbrv}
\bibliography{loewner}
\end{document}